\documentclass[a4paper, 11pt]{article}
\usepackage{amssymb}
\usepackage{float}
\usepackage{amsfonts}
\usepackage{amsthm}
\usepackage{booktabs}
 \usepackage{multirow}
 \usepackage{epstopdf}

\usepackage{amsmath}
\usepackage{graphicx}

\usepackage{hyperref}

\usepackage{subfigure}
\usepackage{authblk}
\usepackage{longtable}
\usepackage{pdflscape}

\newcommand{\RR}{\mathbb{R}}

\newcommand{\PP}{\mathbb{P}}

\newtheorem{thm}{Theorem}
\newtheorem{cor}[thm]{Corollary}
\newtheorem{prop}[thm]{Proposition}

\begin{document}

\title{Chance Constrained Mixed Integer Program: Bilinear and Linear Formulations, and Benders Decomposition}

\author{Bo Zeng,  Yu An and Ludwig Kuznia\\
Dept.\thinspace of Industrial and Management Systems Engineering\\
University of South Florida\\
Tampa, FL 33620\\}
\date{March, 2014}
\maketitle

\abstract{In this paper, we study chance constrained mixed integer
program with consideration of recourse decisions and their incurred cost, developed on a finite
discrete scenario set. Through studying a non-traditional bilinear
mixed integer formulation, we derive its linear counterparts and show that they
could be stronger than existing linear formulations. We also
develop a variant of Jensen's inequality that extends the one for stochastic program.
To solve this challenging problem, we present a variant of Benders
decomposition method in bilinear form, which actually provides an easy-to-use algorithm
framework for further improvements, along with a few
enhancement strategies based on structural properties or Jensen's inequality.
Computational study shows that the presented Benders decomposition method, jointly with appropriate enhancement techniques, outperforms
a commercial solver by an order of magnitude on solving chance constrained program or detecting
its infeasibility. 
}


\section{Introduction}
Chance constrained
mathematical program (CCMP) first
appeared in \cite{charnes1958cost} in 1958 and the (joint) chance
constraint was formally introduced in
\cite{miller.wagner.1965,prekopa1970probabilistic}. It
is often used to capture randomness and to
restrict the associated risk on real system operations \cite{shapiro2009lectures}, and is related to optimization
problems with Value-at-Risk constraints
\cite{sarykalin2008value}.  Mathematically, it can be represented as the following.
Let $\omega$ represents a random
scenario in $\Omega$, $\mathbf{x}$ be the first stage variables and
$\mathbf{y}(\omega)$ be the second stage recourse variables of scenario
$\omega$, we have\begin{align}
\mathbf{CCMP:} \ \ \min\:&\: \mathbf{cx+ F}(\mathbf{y(\omega)})  \label{eq:ccmip} \\
\text{s.t.} \;&\; \mathbf{Ax \geq b} \label{eq:ccmip_firststage} \\
\;&\;  \PP\Big\{ \mathbf{G}(\omega) \mathbf{x} + \mathbf{H}(\omega) \mathbf{y}(\omega) \geq \mathbf{h}(\omega)\Big\} \geq 1-\varepsilon \label{eq:ccmip_chance} \\
        \;&\; \mathbf{x} \in \mathbb{R}^{n1}_+\times\mathbb{Z}^{n2}_+, \mathbf{y(\omega)} \in
        \RR^m_+ \label{eq:ccmip_irv}
\end{align}
where $\mathbf{G}(\omega)$ and $\mathbf{H}(\omega)$
are random matrices,  $\mathbf{h}(\omega)$ is a random column vector of appropriate
dimension, and $\varepsilon$ is the risk tolerance level.
Constraint \eqref{eq:ccmip_chance}, which is denoted as the chance
constraint, enforces that the inequality $\mathbf{G(\omega)} \mathbf{x} +
\mathbf{H}(\omega)
\mathbf{y}(\omega) \geq \mathbf{h}(\omega)$ should be satisfied with a probability
greater than or equal to $1-\varepsilon$. Function $\mathbf{F}$ in the objective
function \eqref{eq:ccmip} represents cost contribution from the
recourse decisions. When the model is static such that cost of
recourse decisions is ignored, function $\mathbf{F}$ is set to $0$
\cite{luedtke2013branch,rusz.pp_w_dd}. Without loss of generality,
we assume that the  set defined by the first stage constraints is not
empty, i.e., $\mathbf{X}=\{\mathbf{x}\in
\mathbb{R}^{n1}_+\times\mathbb{Z}^{n2}_+: \mathbf{Ax} \geq
\mathbf{b}\}\neq \emptyset$.

Because it can explicitly deal with uncertainties and risk
requirements, CCMP has been applied to build
decision making tools for many real systems where uncertainty is a
critical factor in determining system performance. Such applications
include service system design and management
\cite{gurvich2010staffing,noyan2010alternate}, water quality
management \cite{cc.waterqual}, optimal vaccination planning
\cite{tanner2008finding}, energy generation \cite{UC_chance_2004},
and production system design \cite{lejeune2007efficient}. However,
it is well recognized that a CCMP formulation is very challenging to
solve, especially when random $\omega$ is with a general continuous
distribution. In such case, multivariate integration will be
involved in and it is extremely hard to derive a closed-form expression to
represent \eqref{eq:ccmip_chance}, which makes the application
difficult \cite{prekopa2007relationship}. On the other hand, CCMP
with $\Omega$ of a discrete and finite support, i.e.,
$\Omega=\{\omega_1, \dots,\omega_K\}$,  has received very much
attention in recent years. Such CCMP actually
is also constructed  in sampling based methods to solve
those with general continuous distributions
\cite{luedtke2008sample,calafiore2005uncertain}, which makes it
practically applicable. Up to now, several advanced and analytical
solution methods have been developed for this type of CCMP problems.
For the cases where randomness
only appears in right-hand-side and there is no recourse decision,
the concept of $p$-efficient point of the distribution
\cite{prekopa1990dual} is introduced.  It leads to a few methods
that (partially) enumerate those points to derive optimal solutions
\cite{dentcheva2000concavity,beraldi2002branch,lejeune2010mathematical}.
Because the number of efficient points could be huge, such
solution strategy might not be effective in practice. Alternatively,
by using a binary indicator variable and ``big-M'' coefficient, CCMP
can be formulated with a set of linear constraints as the following where $M$ is a
sufficiently large number.
\begin{eqnarray*}
\mathbf{CC-bigM}: \min &&  \mathbf{cx}+ \mathbf{F}(\mathbf{y})\\
\mbox{s.t.} &&  \mathbf{Ax} \geq \mathbf{b}, \\
&& \mathbf{G}_k \mathbf{x} + \mathbf{H}_k \mathbf{y}_k +
Mz_k\geq \mathbf{h}_k, \ k=1,\dots, K\\
&& \sum_{k=1}^K\pi_kz_k\leq \varepsilon\\
&& \mathbf{x} \in \mathbb{R}^{n1}_+\times\mathbb{Z}^{n2}_+,
\mathbf{y}_k \in \RR^m_+, z_k\in\{0,1\}, \ k=1,\dots, K.
\end{eqnarray*}
Note that the enforcement of constraints  of scenario $\omega_k$ is
controlled by $z_k$.  When $z_k=0$,  those constraints must
be satisfied. When $z_k=1$, they can be ignored due to big-M.
Consequently, when $\mathbf{F}(\mathbf{y})$ is $\mathbf{0}$ or linear,
$\mathbf{CC-bigM}$ is a linear  mixed integer program (MIP).
Currently, as it  can handle more general situations, e.g.,
recourse decisions are included and the coefficient matrix is
random, $\mathbf{CC-bigM}$ has become a popular formulation to
address CCMP problems
\cite{rusz.pp_w_dd,luedtke2,luedtke2013branch,tanner2010iis,song2012chance,Shen2013}.
To improve the computational performance on this type of MIP, a
few sophisticated cutting plane methods have been developed
\cite{rusz.pp_w_dd,tanner2010iis,kuccukyavuz2012mixing,luedtke2013branch}.
For example, inequalities based on precedence constrained knapsack
set are derived and included to improve the solution of the MIP
formulation \cite{rusz.pp_w_dd}. Valid inequalities based on IIS
(irreducibly infeasible subsystem) are derived and implemented
within a specialized branch-and-cut procedure \cite{tanner2010iis}.
Based on existing research on the mixing set
\cite{gunluk2001mixing}, strong inequalities generalizing star
inequalities \cite{atamturk2000mixed} are derived to strengthen
the instance only with random right-hand-side
\cite{luedtke2,kuccukyavuz2012mixing}. Such research is further
generalized to compute those with a random coefficient matrix in
\cite{luedtke2013branch}, where strong inequalities that aggregate
Benders feasibility cuts are developed using the mixing set
structure.

The aforementioned algorithms and methods, especially MIP based
approaches, have significantly improve our solution capability to
compute CCMPs. Nevertheless, those algorithms are  complicated for
practitioners to use. Also, they often depend on some non-trivial assumptions and are not applicable for solving
general type instances, e.g., where cost from recourse decisions should
be considered or constraints on recourse decisions could make the
whole formulation infeasible. Indeed, compared to the most relevant
stochastic program (SP), which has received numerous research
efforts and can be solved by various efficient algorithms for
real applications
\cite{shapiro2009lectures,birge2011introduction}, solution methods
to CCMP are rather very limited and less capable. Given the
essential connection between SP and CCMP, we believe that it would be
an ideal strategy to make use of the exiting large amount of research
results on SP to investigate and develop general and effective
algorithms for CCMP, which however receives little attention yet.
Hence, in this paper, we seek to make progress toward this end by
considering solving CCMP with  a set of discrete and finite
scenarios.

We begin our exposition in Section \ref{sect_formulation} by
presenting a non-traditional bilinear formulation of CCMP, strong
linear reformulations, and Jensen's inequality for CCMP. Inspired by
using Benders decomposition method to solve SP, the arguably most
popular approach for SP, we then in Section \ref{sect_Benders}
provide a bilinear Benders reformulation and a bilinear variant of
Benders decomposition method. We also design a few enhancement strategies based on structural 
properties or Jensen's inequality, all
of which can be easily implemented. In Section \ref{Sect_numerical},
we perform a numerical study on a set of unstructured random instances and 
a set of operating room scheduling instances  arising from healthcare applications. Our results 
show that the presented bilinear Benders decomposition method, jointly with appropriate enhancement 
strategies,   
typically
outperforms the state-of-the-art commercial solver by an order of
magnitude on deriving optimal solutions or reporting infeasibility.
Finally, we provide concluding
remarks in Section \ref{Sect_conclusion}.


We mention that our paper provides the following contributions and
insights, some of which might be counterintuitive. $(i)$ The
bilinear formulation of CCMP is informative. Its linearized
formulations can be theoretically stronger and computationally more
friendly than the popular $\mathbf{CC-bigM}$ formulation. So, it is worth
more research efforts to investigate the mathematical structure and
to study fast computing algorithms. $(ii)$ A variant
of Jensen's inequality is derived for CCMP, which generalizes its SP
counterpart. We demonstrate that it could also generate a significant computational improvement on 
solving difficult real problems. $(iii)$ Different from existing understanding that
Benders decomposition may not be a good method \cite{tanner2010iis},
Benders decomposition, in the presented bilinear form, is very
capable to solve CCMP. Indeed, it is interesting to note that
Benders decomposition could be more effective (in terms of
iterations) in solving CCMP than in solving SP. $(iv)$ A few
non-trivial enhancement strategies that specifically consider the
structure of CCMP are developed, which could lead to further
performance improvement. $(v)$ The presented Benders decomposition
method yields an easy-to-use fast algorithm strategy that basically does not depend on special
assumptions and can solve
general CCMPs, which have not been addressed in existing
literature. Note that it rather provides a framework that can
incorporate numerous existing results on Benders decomposition for
SP to improve our solution capacity of CCMP. 

\section{CCMP: Bilinear Formulations and Properties}
\label{sect_formulation}
\subsection{A  Bilinear Formulation of CCMP}
\label{sect_formulation_assum} To simplify our exposition, we define
a scenario contributing to the satisfaction of the chance constraint
as a responsive scenario,
and a scenario that can be ignored (i.e., the corresponding
constraints in that scenario can be violated) as a non-responsive
scenario. In this paper,
we assume that the (expected) recourse cost is derived as the
weighted sum of costs from responsive scenarios, which is
 $\mathbf{F(y)}= \sum_{k=1}^K
\pi_k(\mathbf{f}_k\mathbf{y}_k)(1-z_k)$, where $\pi_k$ is the
realization probability and $\mathbf{f}_k$ is cost vector of
$\mathbf{y}_k$, for scenario $\omega_k$, $k=1,\dots, K$. Next, different
from the popular $\mathbf{CC-bigM}$ MIP formulation, we present a mixed integer
bilinear formulation of CCMP as the following.
\begin{align}
\mathbf{CC-MIBP:} \ \theta^*(\varepsilon) = \quad \min \;&\; \mathbf{cx} +
\sum_{k=1}^K \pi_k \eta_k \label{mip_cc_mibp}\\
\text{s.t.}  \;&\; \mathbf{Ax} \geq \mathbf{b}   \label{mibp_eq_first} \\
\;&\;  \eta_k = \mathbf{f}_k\mathbf{y}_k(1-z_k), \ k=1,\dots,K \label{eq_sub_obj_pri}\\
  \;&\; (\mathbf{G}_k \mathbf{x} + \mathbf{H}_k\mathbf{y}_k - \mathbf{h}_k)(1-z_k)\geq 0, \ k = 1,\dots,K \label{mibp_scen_const}\\
\;&\; \sum_{k=1}^K \pi_k z_k \leq \varepsilon   \label{mibp_chance_const}\\
\;&\; \mathbf{x} \in \mathbb{R}^{n1}_+\times\mathbb{Z}^{n2}_+, \
\mathbf{y}_k \in \RR^m_+, \;\; z_k \in \{0,1\}, \ k=1,\dots,K.
\label{mibp_mip_cc_end}
\end{align}
The validity of  $\mathbf{CC-MIBP}$ is obvious. Note that
by assigning $z_{k}$ to one or zero, the impact of scenario $k$, including the
cost contribution of recourse decisions and feasibility requirements
from recourse constraints,  will be removed from the whole
formulation or imposed explicitly. Although  $\mathbf{CC-MIBP}$ is a non-traditional nonlinear
formulation,
it provides a  more direct
and natural representation of the underlying combinatorial structure, which captures CCMP without information loss
and allows us to gain deep
insights.

We define
$\theta^*(\varepsilon)=+\infty$ if it is infeasible and
$\theta^*(\varepsilon)=-\infty$ if it  is unbounded. Noting that an
instance  with a larger $\varepsilon$ is a relaxation of another one
with a smaller $\varepsilon$, the next result follows.
\begin{prop}
\label{prop_monotone}
For $0\leq
\varepsilon_1\leq \varepsilon_2\leq 1$, we have
$\theta^*(\varepsilon_1)\geq \theta^*(\varepsilon_2)$.
\end{prop}
We mention that when $\varepsilon=0$, it is clear that $z_k=0$ for all
$k$ and $\mathbf{CC-MIBP}$ reduces to a linear MIP, which actually
is the underlying SP model. According to Proposition
\ref{prop_monotone}, the optimal value of SP, i.e., $\theta^*(0)$,
provides an upper bound to $\theta^*(\varepsilon)$ for any
$\varepsilon>0$.\\



\noindent \textbf{Remarks:} We note that it is often assumed in the
study of SP that it has a finite optimal value, which could help
researchers focus on solving real life problems that are generally
feasible and bounded. Indeed, given the fact that SP is a linear
program or mixed integer program, which has many advanced
preprocessing techniques developed to detect infeasibility or
unboundedness, such assumption might not be restrictive. However,
the situation of $\mathbf{CC-MIBP}$ (or CCMP in general) is more involved. Different
$\varepsilon$ may cause $\mathbf{CC-MIBP}$ infeasible, finitely optimal, or
unbounded. For example, when $\varepsilon=0$, the resulting SP model
must deal with constraints in all scenarios, which could conflict to
each other, forcing the whole formulation infeasible. When
$\varepsilon$ increases to a larger value, we can control binary
$\mathbf{z}$ variables to remove conflicting scenarios so that
constraints from remaining scenarios  are compatible, rendering the
formulation feasible.
\textrm{}\\


\noindent \textbf{Assumption:} In this paper, to help us focus on
typical problems, we make an assumption that is independent of
$\varepsilon$. We think such independence is of a critical
value for real system users as the sensitivity analysis, including
feasibility and the total cost, with respect to $\varepsilon$
reveals a fundamental feature of CCMP on managing risks and
uncertainty. Specifically, we assume that for any scenario
$\omega_k$, there $\exists \ \mathbf{x}^0\in \mathbf{X}$ such that the
recourse problem is finitely optimal, i.e.,
\begin{eqnarray}
\label{eq_assumption} -\infty
<\eta_k^*(\mathbf{x}^0)=\min\{\mathbf{f}_k\mathbf{y}_k:
\mathbf{H}_k\mathbf{y}_k \geq \mathbf{h}_k - \mathbf{G}_k
\mathbf{x}^0, \mathbf{y}_k \in \RR^m_+\}<+\infty.
\end{eqnarray}
If such assumption is violated in
$\omega_k$, there exist three cases which can be processed by the following operations. \\
\emph{Case $(i)$}: $\eta^*_k(\mathbf{x}) =+\infty, \forall \
\mathbf{x} \in\mathbf{X}$. In this case, we fix $z_k=1$ in
$\mathbf{CC-MIBP}$. Alternatively, we can update
$\varepsilon=\varepsilon-\pi_k$ and remove constraints and recourse
variables of $\omega_k$ from $\mathbf{CC-MIBP}$. \\
\emph{Case
$(ii)$}: $\eta^*_k(\mathbf{x})=-\infty, \forall \ \mathbf{x}
\in\mathbf{X}$. In this case, we fix $z_k=0$ in $\mathbf{CC-MIBP}$.
Alternatively, we can simply remove constraints and variables of
$\omega_k$ from $\mathbf{CC-MIBP}$. If $\mathbf{CC-MIBP}$ with
remaining scenarios is feasible, the original formulation is
unbounded. Otherwise, both
the updated one and original one are infeasible. \\
\emph{Case $(iii)$}: $\mathbf{X}=\mathbf{X}_1\cup\mathbf{X}_2$ such
that $\eta^*_k(\mathbf{x})=+\infty, \forall \ \mathbf{x}
\in\mathbf{X}_1$ and $\eta^*_k(\mathbf{x})=-\infty, \forall \ \mathbf{x} \in\mathbf{X}_2$. Basically, such a case is very rare as
$\eta_k$ does not have a transition between $\mathbf{X}_1$ and
$\mathbf{X}_2$. Indeed, when $\mathbf{X}$ is  a feasible set of LP,
this situation will not occur. If it does happen, we make use of
Branch-and-Bound technique by creating two branches for
$\mathbf{X}_1$ and $\mathbf{X}_2$ respectively. Then, operations
presented for Case $(i)$ and Case $(ii)$ can be applied to those
branches for simplification.

Hence,  we can intuitively interpret the aforementioned assumption as
that it ensures all scenarios behave in a regular fashion because there
is no scenario simply causing unboundness or becoming irrevelent.
Cases that do not comply with that assumption are rather extreme and
they  can be dealt with in the preprocess stage. So, our
assumption imposes little restriction on CCMP in practice.


\subsection{Deriving Linear Reformulations from Bilinear Formulation}
$\mathbf{CC-MIBP}$'s original bilinear form might not be
computationally friendly. It would be ideal to convert it into a
linear formulation. In particular, we are interested in a tight
linear formulation stronger than the popular $\mathbf{CC-bigM}$.

By expanding \eqref{mibp_scen_const}, we have
\begin{eqnarray}
\label{mibp_scen_const_1} \mathbf{G}_k
\mathbf{x}-\mathbf{G}_k\mathbf{x}z_k + \mathbf{H}_k\mathbf{y}_k-
\mathbf{H}_k\mathbf{y}_kz_k + \mathbf{h}_kz_k\geq \mathbf{h}_k, \
k=1,\dots,K.
\end{eqnarray}
Consider a situation where no recourse decision is involved, i.e.,
$m=0$ in \eqref{mibp_mip_cc_end} and \eqref{eq_sub_obj_pri}
disappears, and $\mathbf{G}_k=\mathbf{G}$ for all $k$, i.e., only
right-hand-side $\mathbf{h}_k$ is random. According to
\cite{kuccukyavuz2012mixing}, we can assume without loss of
generality that $\mathbf{h}_k\geq \mathbf{0}$ for $k=1,\dots, K$.
Note that, no matter which scenarios are responsive in determining
an optimal solution, we have $\mathbf{Gx}\geq
\min_{l=1,\dots,K}\{\mathbf{h}_l\}$, where $\min$ is applied in a component-wise fashion.
Hence, projecting out $\mathbf{Gx}z_k$ through replacing it by
$\min_{l=1,\dots,K}\{\mathbf{h}_l\}z_k$ in \eqref{mibp_scen_const_1}, we have
\begin{eqnarray}
\label{mibp_scen_const_2} \mathbf{Gx}-\min_{l=1,\dots,K}\{\mathbf{h}_l\}z_k +
\mathbf{h}_kz_k\geq \mathbf{h}_k, \ \ k=1,\dots, K.
\end{eqnarray}
Because $\mathbf{Gx}\geq \min_k\{\mathbf{h}_k\}\geq \mathbf{0}$, it
simply leads to the next result.
\begin{prop}
\label{prop_dominant} The linear model
\begin{align}
\label{eq_bilinear_linear1} \theta^*(\varepsilon) = \min
\Big\{\mathbf{cx}: \eqref{mibp_eq_first}, \eqref{mibp_chance_const},
\eqref{mibp_mip_cc_end}, \ \mathbf{Gx} +
(\mathbf{h}_k-\min_{l=1,\dots,K}\{\mathbf{h}_l\})z_k\geq \mathbf{h}_k, \ \forall k \Big\}
\end{align}
is a valid formulation for CCMP (without recourse opportunities).
Moreover, it dominates the following linear formulation
\begin{align}
\theta^*(\varepsilon) = \min \{\mathbf{cx}: \eqref{mibp_eq_first},
\eqref{mibp_chance_const}, \eqref{mibp_mip_cc_end}, \ \mathbf{Gx}+
\mathbf{h}_kz_k\geq \mathbf{h}_k, \ \forall k \},
\end{align}
which is adopted as CCMP formulations in
\cite{luedtke2,kuccukyavuz2012mixing}. \hfill $\square$
\end{prop}

We can further extend results in Proposition \ref{prop_dominant} to
more general cases where CCMP has recourse decisions. Assume that
the set defined by \eqref{mibp_scen_const_1}   can be equivalently
represented by the following linear inequality:
\begin{eqnarray}
\label{eq_linear_dummy} \mathbf{G}_k \mathbf{x}+
\mathbf{H}_k\mathbf{y}_k +\mathbf{h}_kz_k\geq
\mathbf{h}_k+\mathbf{q}_kz_k
\end{eqnarray}
 where $\mathbf{q}_k$ is
a coefficient vector with an appropriate dimension. Obviously, to
ensure the equivalence between \eqref{mibp_scen_const_1} (or
\eqref{mibp_scen_const}) and \eqref{eq_linear_dummy}, we should have
for any $\mathbf{x}\in \mathbf{X}$, there exists $\mathbf{y}_k$ such
that the inequality $\mathbf{G}_k \mathbf{x}+
\mathbf{H}_k\mathbf{y}_k\geq \mathbf{q}_k$ is satisfied when $z_k=1$
(noting that the cost contribution of $\mathbf{y}_k$ is not a
concern when $z_k=1$). Otherwise, additional restriction will be
imposed on $\mathbf{X}$, which suggests that the impact of scenario
$\omega_k$ is not removed when $z_k=1$. Hence, $q_{k,i}$, the
$i^{\textrm{th}}$ component of $\mathbf{q}_k$, should have
\begin{eqnarray}
\label{eq_bilevel} q_{k,i}\leq \min_{\mathbf{x\in \mathbf{X}}}
\max_{\mathbf{y}_k\in \mathbf{Y}_k} \
 (\mathbf{G}_k \mathbf{x}+
\mathbf{H}_k\mathbf{y}_k)_i,
\end{eqnarray}
for every $i$, where $\mathbf{Y}_k$ represents the set defined by
variable bounds of $\mathbf{y}_k$.
 Certainly, setting
$q_{k,i}$ to achieve the equality will give us the tightest linear
inequality to replace \eqref{mibp_scen_const}. Nevertheless, it
involves computing a bilevel robust optimization problem, which is
computationally expensive. A simpler strategy is to fix
$\mathbf{y}_k$ to $\mathbf{y}_k^0\in \mathbf{Y}_k$ and let $q_{k,i}=
\min\{(\mathbf{G}_k\mathbf{x})_i: \mathbf{x}\in
\mathbf{X}\}+\{\mathbf{H}_k\mathbf{y}^0_k\}_i$. Actually, for a situation
where $\mathbf{0}\in \mathbf{Y}_k$ and $\mathbf{f}_k\geq
\mathbf{0}$, a strong model only with linear constraints could be
derived to represent CCMP.
\begin{prop}
\label{prop_linear_bilinear} When $\mathbf{0}\in \mathbf{Y}_k$, let
$q^*_{k,i}= \min\{(\mathbf{G}_k\mathbf{x})_i: \mathbf{x}\in
\mathbf{X}\}$ for all $i$. If $\mathbf{q}^*_k>-\infty$ and
$\mathbf{f}_k\geq \mathbf{0}$ for all $k$, the following linear
model is a valid formulation for CCMP:
\begin{align}
\label{eq_bilinear_linear2} \theta^*(\varepsilon) = \min
\{\mathbf{cx} + \sum_{k=1}^K \pi_k \mathbf{f}_k\mathbf{y}_k:
\eqref{mibp_eq_first}, \eqref{mibp_chance_const},
\eqref{mibp_mip_cc_end}, \ \mathbf{G}_k \mathbf{x} +
\mathbf{H}_k\mathbf{y}_k + (\mathbf{h}_k- \mathbf{q}^*_k) z_k \geq
\mathbf{h}_k, \ \forall k\}
\end{align}
\end{prop}
\begin{proof}

For any $\mathbf{x}$, with $\mathbf{f}_k\geq \mathbf{0}$, i.e., recourse decisions incur
non-negative cost, it would be optimal to set
$\mathbf{y}_k=\mathbf{0}$ assuming it is feasible. Under such a
situation, it is valid to simplify \eqref{eq_sub_obj_pri} to
$\eta_k=\mathbf{f}_k\mathbf{y}_k$, which ensures there is no cost
contribution when $z_k=1$. 

Indeed, according to the definition of $\mathbf{q}^*_k$, we have
$\mathbf{q}^*_k - \mathbf{G}_k \mathbf{x}\leq \mathbf{0}$ for any
$\mathbf{x} \in \mathbf{X}$. Hence, given that $\mathbf{0}\in
\mathbf{Y}_k$, it is guaranteed that $\mathbf{0}$ is a feasible
point of $\{\mathbf{y}_k \in \mathbf{Y}_k: \mathbf{H}_k\mathbf{y}_k
\geq \mathbf{q}^*_k - \mathbf{G}_k \mathbf{x}\}$ for any $\mathbf{x}
\in \mathbf{X}$. Therefore, \eqref{eq_bilinear_linear2} is a valid
linear formulation of CCMP.
\end{proof}
\noindent\textbf{Remarks:}\\
$(i)$ Clearly, applying the aforementioned results will produce a
strong linear formulation of CCMP, which alleviates our concern on
big-M coefficients. Note that deriving $q^*_{k,i}$ requires
computing an MIP for every $k$ and $i$, which may incur non-trivial
computational expense. One strategy is to relax those MIPs as LPs and derive weaker coefficients with less computational
expense. To compute large-scale practical instances with many
scenarios and constraints, we believe that computationally more
friendly methods are necessary to achieve a balance between the
coefficient quality
and computational time. In addition, it is worth mentioning that those results
do not depend on the variable types of $\mathbf{y}_k$ for any $k$. Hence, it can be used
to derive strong linear formulations for those whose recourse problems could be mixed integer programs.
\\
$(ii)$ Our discussions also show that for a case where
$\mathbf{q}^*_{k,i}$ is not finite for some $k$ and $i$ or cost of
some recourse decision is not non-negative, it will be more
complicated. We might not be able to derive a finite coefficient for
$z_k$ in \eqref{eq_bilinear_linear2}, which also implies big-M
method might not work under those situations. Hence, a deeper study
on those cases is needed to build strong formulations. On
this line, an analytical study on deriving strong formulations and
valid inequalities is presented in \cite{MingBoKai}.

Instead of using computational methods to derive linear
formulations, a simpler and more general approach is to directly
linearize $\mathbf{CC-MIBP}$ and obtain an MIP in a higher
dimension. Given the fact that $z_k$ is binary,
constraints in \eqref{eq_sub_obj_pri} and \eqref{mibp_scen_const_1}
(equivalently \eqref{mibp_scen_const})
  can be easily
converted into linear constraints by applying \emph{McCormick
linearization method} on bilinear terms $\mathbf{x}z_k$ and
$\mathbf{y}_kz_k$ \cite{mccormick1976computability}. For example,
say $x_j$ is one component of $\mathbf{x}$. Its bilinear term
$x_jz_k (=\lambda_{jk})$ can be linearized by including the following constraints
\begin{eqnarray}
\label{eq_linearizationmethod} \lambda_{jk}\leq x_j, \ \lambda_{jk}\leq
U_{x_j}z_k, \ \lambda_{jk}\geq x_j+U_{x_j}(z_k-1), \ \lambda_{jk} \geq 0,
\end{eqnarray} where $U_{x_j}$ is an upper bound of single variable
$x_j$. We recognize that such linearization method depends on an
upper bound of each individual variable. Nevertheless, it is
typically much easier to estimate a strong  upper bound for an
individual variable, which indicates the linear formulation obtained
using McCormick linearization method will be tighter than $\mathbf{CC-bigM}$. Indeed, in
many practices, upper bounds of individual variables (e.g., binary
or bounded variables) are naturally defined or available
\cite{gurvich2010staffing,tanner2010iis}. In such a case, there is
no need to introduce a very big number as the variable upper bound
when we perform linearization. Furthermore, this linearization
strategy can deal with situations where recourse decisions are with
unrestricted cost coefficients, e.g., recourse decisions are with
negative cost. So, because of its generality and simplicity, in the
remainder of this paper, unless explicitly mentioned,
$\mathbf{CC-MIBP}$ and its linearized counterpart obtained through
McCormick linearization method are interchangeably used.

Given the rich and original  structural information implied in
$\mathbf{CC-MIBP}$, we believe that more mathematical analysis
should be done on this bilinear formulation and specialized
algorithms need to be developed. Next, following this line, we
present a generalization  of Jensen's inequality in CCMP and
 a variant of Benders decomposition  customized according to this bilinear
form.

\subsection{Jensen's inequality for CCMP}
\label{sect_Jensen} The combinatorial structure explicitly presented
in (\ref{eq_sub_obj_pri}-\ref{mibp_scen_const}) helps us to
extend the classical Jensen's inequality
\cite{jensen1906fonctions,madansky1960inequalities}, a strong
inequality for stochastic programs, to CCMP. Consider the situation
where
 $\mathbf{H}_k=\mathbf{H}$ and
$\mathbf{f}_k=\mathbf{f}$ for all $k$.

\begin{thm}
\label{thm_Jensens} (Extended Jensen's inequality for CCMP) Let
$\mathbf{z}^0$ be an optimal $\mathbf{z}$ with respect to
$\mathbf{x}^0\in \mathbf{X}$ and define $E\eta(\mathbf{x}^0)=
\sum_{k=1}^K\pi_k(1-z^0_k)\eta^*_k(\mathbf{x}^0)$, i.e., the optimal
expected recourse cost from responsive scenarios. If
$E\eta(\mathbf{x}^0)<+\infty$, we have
\begin{eqnarray}
 E\eta(\mathbf{x}^0) \geq
\min\Big\{(1-\sum_{k=1}^K\pi_kz_k) \mathbf{f}\mathbf{\overline y}: &&
\mathbf{H}\mathbf{\overline y}\geq
\frac{\sum_{k=1}^K\pi_k(1-z_k)(\mathbf{h}_k-\mathbf{G}_k\mathbf{x}^0)}{1-\sum_{k=1}^K\pi_kz_k},
\sum_{k=1}^K \pi_k z_k \leq \varepsilon, \notag\\
&& \mathbf{\overline y}\geq
\mathbf{0}, \mathbf{z}\in \{0,1\}^K\Big\}. \label{eq_jensen_chance}
\end{eqnarray}
\end{thm}
\begin{proof}
We define $\tilde \pi_k =
\frac{\pi_k(1-z^0_k)}{1-\sum_{k=1}^K\pi_kz^0_k}$. Note that $\tilde
\pi_k\geq 0$ and
$$\sum_{k=1}^K \tilde \pi_k= \sum_{k=1}^K \frac{\pi_k(1-z^0_k)}{1-\sum_{k=1}^K\pi_kz^0_k}
=\frac{\sum_{k=1}^K \pi_k -\sum_{k=1}^K \pi_kz^0_k}{1-\sum_{k=1}^K\pi_kz^0_k}=1.$$

Furthermore, we define $\tilde \eta(\omega) = \min\{\mathbf{fy}:
\mathbf{Hy}\geq \mathbf{h}(\omega)-\mathbf{G}(\omega)\mathbf{x}^0,
\mathbf{y}\geq \mathbf{0}\}$, where $\omega\in
\{\omega_1,\dots,\omega_K\}$. According to the duality
theory, $\tilde \eta(\omega) = \max\{
(\mathbf{h}(\omega)-\mathbf{G}(\omega)\mathbf{x}^0)^T\mathbf{u}:
\mathbf{H}^T\mathbf{u}\leq \mathbf{f}, \mathbf{u}\geq \mathbf{0}\}.$
Because $\mathbf{H}$ and $\mathbf{f}$ are independent of scenarios,
we have that $\tilde \eta (\omega)$ is a convex function of random
$\omega$. Therefore, following the definition of $\tilde \pi_k$, we
have
$$\sum_{k}\tilde \pi_k\tilde \eta_k \geq  \min\Big\{\mathbf{f \overline y}:
\mathbf{H \overline y}\geq \sum_k\tilde \pi_k(\mathbf{h}_k-\mathbf{G}_k\mathbf{x}^0),
\mathbf{\overline y}\geq \mathbf{0}\Big\}. $$
Multiplying both sides by $(1-\sum_{k=1}^K\pi_kz^0_k)$, we have
$$E\eta(\mathbf{x}^0)=\sum_{k}\pi_k(1-z^0_k)\tilde \eta_k \geq
\min\Big\{ (1-\sum_{k=1}^K\pi_kz^0_k)\mathbf{f \overline y}: \mathbf{H
\overline y}\geq \sum_k\tilde
\pi_k(\mathbf{h}_k-\mathbf{G}_k\mathbf{x}^0), \mathbf{\overline
y}\geq \mathbf{0}\Big\}.$$ Also, noting that $\mathbf{z}^0$ is a
particular $\mathbf{z}$ satisfying $\sum_{k=1}^K \pi_k z_k \leq
\varepsilon$, we have
\begin{eqnarray*}
 &&\min \Big\{(1-\sum_{k=1}^K\pi_kz^0_k) \mathbf{f \overline y}:
\mathbf{H \overline y}\geq \sum_k\tilde \pi_k(\mathbf{h}_k-\mathbf{G}_k\mathbf{x}^0),
\mathbf{\overline y}\geq \mathbf{0}\Big\} \\
&&  \geq \min \Big\{ (1-\sum_{k=1}^K\pi_kz_k)  \mathbf{f}\mathbf{\overline y}:
\mathbf{H}\mathbf{\overline y}  \geq \frac{\sum_{k=1}^K\pi_k(1-z_k)(\mathbf{h}_k-\mathbf{G}_k\mathbf{x}^0)}{1-\sum_{k=1}^K\pi_kz_k},
\sum_{k=1}^K \pi_k z_k \leq \varepsilon, \\
&& \ \ \quad \quad \mathbf{\overline y}\geq
\mathbf{0}, \mathbf{z}\in \{0,1\}^K\Big\}.
\end{eqnarray*}
Therefore, the desired inequality follows.
\end{proof}

\noindent \textbf{Remarks:}\\
When $\varepsilon=0$, we have $z_k=0$ for all $k$ in both
$\mathbf{CC-MIBP}$ and \eqref{eq_jensen_chance}.  Note that under
such situation, CC-MIBP is SP formulation and
\eqref{eq_jensen_chance} reduces to
\begin{eqnarray}
\label{eq_jensen_SP} E\eta(\mathbf{x}^0) \geq
\min\Big\{\mathbf{f}\mathbf{\overline y}:
\mathbf{H}\mathbf{\overline y}\geq
\sum_{k=1}^K\pi_k(\mathbf{h}_k-\mathbf{G}_k\mathbf{x}^0),
\mathbf{\overline y}\geq \mathbf{0}\Big\},
\end{eqnarray}
which  is
Jensen's inequality for SP. Nevertheless, as it involves solving an
MIP with fractional inequalities, \eqref{eq_jensen_chance} in
general is much more complicated than its SP counterpart.  Those
fractional inequalities can be converted into equivalent bilinear
inequalities as follows
\begin{eqnarray}
\label{eq_jensen_chance_bilinear}
\mathbf{H}\mathbf{\overline y}(1-\sum_{k=1}^K\pi_kz_k)\geq
\sum_{k=1}^K\pi_k(1-z_k)(\mathbf{h}_k-\mathbf{G}_k\mathbf{x}^0).
\end{eqnarray}
Hence, with fractional inequalities replaced by \eqref{eq_jensen_chance_bilinear},
the minimization problem in the right-hand-side of \eqref{eq_jensen_chance}
becomes an MIP with bilinear inequalities.
Actually, for some special case, those fractional or bilinear inequalities reduce
to linear ones.
\begin{cor}
\label{cor_jensens} Assume $\pi_k=\frac{1}K$ for all $k$ and
$\mathbf{f}\geq \mathbf{0}$. Let $L$ be the integer such that
$\frac{L}{K}\leq \varepsilon < \frac{L+1}{K}$. The extended Jensen's
inequality in \eqref{eq_jensen_chance} can be simplified as
\begin{eqnarray}
\label{eq_linear_jensen} E\eta(\mathbf{x}^0) \geq
\min\Big\{\frac{K-L}{K}\mathbf{f\overline y}:
\mathbf{H}\mathbf{\overline y}\geq
\sum_{k=1}^K\frac{1-z_k}{K-L}(\mathbf{h}_k-\mathbf{G_kx}^0),
\sum_{k=1}^K z_k = L, \mathbf{\overline y}\geq \mathbf{0},
\mathbf{z}\in \{0,1\}^K\Big\}.
\end{eqnarray}
\end{cor}
\begin{proof}
When $\pi_k=\frac{1}K$ for all $k$ and $\mathbf{f}\geq \mathbf{0}$,
it is easy to see that there is an optimal $\mathbf{z}^0$ with
$\sum_{k}z^0_k=L$ for  given $\mathbf{x}^0$. Hence, based on the
proof of Theorem \ref{thm_Jensens}, we replace $\pi_k$ by $\frac 1K$
and replace $\sum_{k}z^0_k$ by $L$ in  \eqref{eq_jensen_chance}.
Then, the inequality in \eqref{eq_linear_jensen} follows.
\end{proof}

We mention that Jensen's inequality is not only of a theoretical
contribution, but also could bring significantly computational
benefits. Through replacing $\mathbf{x}^0$ by variable $\mathbf{x}$,
the inequality in \eqref{eq_jensen_chance} or
\eqref{eq_linear_jensen} can be converted into a valid inequality to
bound the optimal expected recourse cost from below. Application of
such strategy in SP formulations has been proven to be very
successful in computing complicated real problems
\cite{erdogan2013dynamic,batun2011operating}. Hence, it would be
interesting to investigate its performance in solving CCMP problems.

\noindent \textbf{Remarks:} \\
One may think that we can derive a variant of Jensen's inequality by
treating $\mathbf{CC-bigM}$ as an SP problem with both $\mathbf{x}$
and $\mathbf{z}$ being first stage variables. Nevertheless, such
 strategy generally leads to  a trivial inequality that is of little interest. We
 illustrate our argument by considering a given first stage solution
 $(\mathbf{x}^0,\mathbf{z}^0)$ under the situation described in
 Corollary \ref{cor_jensens}.

 According to \eqref{eq_jensen_SP} and formulation $\mathbf{CC-bigM}$, we have
 \begin{eqnarray*}
 E\eta(\mathbf{x}^0) & \geq
& \min\Big\{\mathbf{f}\mathbf{\overline y}:
\mathbf{H}\mathbf{\overline y}\geq
\frac{1}{K}\sum_{k=1}^K(\mathbf{h}_k-\mathbf{G}_k\mathbf{x}^0-Mz^0_k),
\mathbf{\overline y}\geq \mathbf{0}\Big\}\\
&=& \min\Big\{\mathbf{f}\mathbf{\overline y}:
\mathbf{H}\mathbf{\overline y}\geq
\frac{1}{K}\sum_{k=1}^K(\mathbf{h}_k-\mathbf{G}_k\mathbf{x}^0)-\frac{M}K\sum_{k=1}^Kz^0_k,
\mathbf{\overline y}\geq \mathbf{0}\Big\}\\
&=& \min\Big\{\mathbf{f}\mathbf{\overline y}:
\mathbf{H}\mathbf{\overline y}\geq
\frac{1}{K}\sum_{k=1}^K(\mathbf{h}_k-\mathbf{G}_k\mathbf{x}^0)-\frac{ML}{K},
\mathbf{\overline y}\geq \mathbf{0}\Big\}\\
&=& 0.
 \end{eqnarray*}
 The second to the last equality follows the fact that $\sum_{k}z^0_k=L$ in one
 optimal solution. Furthermore, when $M$ is sufficiently
 large and $L>0$, constraints inside the minimization problem are always satisfied
 by any $\mathbf{\overline y}\geq \mathbf{0}$. Given that $\mathbf{f}\geq 0$,
 an obvious optimal solution is $\mathbf{y}^*=0$. Hence,
  the last equality follows, which is trivial. Again, this
 illustration shows that the linear formulation with big-M is weak
 in deriving structural properties.

%
%
%

\section{Benders Decomposition for $\mathbf{CC-MIBP}$}
\label{sect_Benders}
\subsection{Bilinear Benders Reformulation and Decomposition Method}
\label{sect_Benders_refor}
Given the definition of $z_k$ and bilinear formulation
$\mathbf{CC-MIBP}$ defined in
(\ref{mip_cc_mibp}-\ref{mibp_mip_cc_end}), we can naturally present
its Benders reformulation in a bilinear form. Consider the dual
problem of recourse problem in scenario $\omega_k$ as follows.
\begin{eqnarray}
\mathbf{DP}_k: \  \ \quad \max && (\mathbf{h}_k-\mathbf{G}_k\mathbf{x}^0)^T \mathbf{u}_k \label{eq_sub_obj}\\
\mbox{s.t.} && \mathbf{H}^T_k \mathbf{u}_{k}\leq \mathbf{f}_k \\
&& \mathbf{u}_k\in \mathbb{R}'_+. \label{eq_sub_var}
\end{eqnarray}
According to our assumption made in Section
\ref{sect_formulation_assum} and the duality theory, it is clear
that the feasible set of this dual problem is non-empty, regardless
of $\mathbf{x}^0$. Denote all the extreme points of this set  by
$E^k=\{\mu^1,\dots, \mu^{p_k}\}$, all extreme rays by $F^k=
\{\upsilon^1,\dots, \upsilon^{t_k}\}$, and let $E$ and $F$ be the
corresponding collections over all $k$.  Next, we present Benders
reformulation of $\mathbf{CC-MIBP}$.
\begin{thm}
\label{thm_BD_reformulation} (Bilinear Benders Reformulation) The
$\mathbf{CC-MIBP}$ formulation can be equivalently reformulated as:
\begin{align}
\mathbf{BD-MIBP}: \theta^* = \min \;&\; \mathbf{cx} + \sum_k \pi_k\eta_k & \label{mip_cc_BD}\\
\mbox{s.t.}  \;&\; \mathbf{Ax} \geq \mathbf{b} \\
\;&\; (\mathbf{h}_k-\mathbf{G}_k\mathbf{x})^T\mu^l(1-z_k)\leq \eta_k; \
 \mu^l\in E^k, \ k=1,\dots, K & \label{eq_xpoints} \\  
\;&\; (\mathbf{h}_k-\mathbf{G}_k\mathbf{x})^T\upsilon^l(1-z_k) \leq
0; \ \upsilon^l \in F^k, \
k=1,\dots, K & \label{eq_xrays}\\ 
\;&\; \sum_{k=1}^K \pi_k z_k \leq \varepsilon & \label{chance_const}\\
\;&\; \mathbf{x} \in \mathbb{R}^{n_1}_+\times\mathbb{Z}^{n_2}_+ ;  \
z_k \in \{0,1\}, k=1,\dots,K. \label{mip_cc_end}
\end{align}
$(i)$ If all recourse problems are   feasible for any
$\mathbf{x}\in \mathbf{X}$, only those in \eqref{eq_xpoints}, i.e.,
those defined by extreme points, are needed in $\mathbf{BD-MIBP}$;
$(ii)$ if recourse decisions are with zero cost, only constraints in
\eqref{eq_xrays}, i.e., those defined by extreme rays, are needed in
$\mathbf{BD-MIBP}$.
 \hfill $\square$
\end{thm}
Again, bilinear terms in constraints in \eqref{eq_xpoints} and
\eqref{eq_xrays} can be linearized using McCormick linearization
method, which converts  \eqref{eq_xpoints} and \eqref{eq_xrays} into
linear inequalities and the whole formulation into an MIP. We
mention that, comparing to the original $\mathbf{CC-MIBP}$
formulation, this bilinear Benders reformulation
(and its linearized counterpart) has a couple of advantages. \\

\noindent \textbf{Remarks:} \\
\noindent $(i)$ We observe that only the first stage decision
variables $\mathbf{x}$ are involved in this bilinear Benders
reformulation and there is no recourse variable. Such an observation
actually gives us an essential advantage in eliminating the negative
impact of big-M parameters in linearization. As $\mathbf{x}$ is
typically introduced to model decisions with concrete definitions in
practice, e.g., quantity of available resources or binary decisions,
their upper bounds can often be obtained from specific applications.
Actually, we note that it is straightforward to obtain physically
valid upper bounds for the first stage decision variables in the
majority of existing CCMP formulations adopted for practical
problems. See
\cite{UC_chance_2004,gurvich2010staffing,noyan2010alternate,shen2010expectation}
for a few examples. Therefore, there is no ``big-M'' issue and our
$\mathbf{BD-MIBP}$ formulation is strong, which might lead to
computational improvements in solving real problems.\\
\noindent $(ii)$ Note also from \eqref{eq_xpoints} and
\eqref{eq_xrays} (and their linearized counterparts) that $z_k$ acts
on all constraints generated from extreme points and rays of
subproblem/scenario $k$. Its binary property indicates the inclusion
or omission of a whole group of constraints. If $\mathbf{z}$
variables are fixed, the resulting problem has  a structure same as
the Benders reformulation of a stochastic program. Hence, it
suggests that,  if  Benders decomposition method works well for the
underlying SP formulation, it should work well for CCMP. 

The next result follows directly from Theorem
\ref{thm_BD_reformulation} and it provides the basis to design our
Benders decomposition method.
\begin{prop}
\label{prop_validity_Benders} Let $\hat E^k$ and $\hat F^k$ be
subsets of extreme points and extreme rays of $\mathbf{DP}_k$
for all $k$, and $\hat E$ and $\hat F$ be the corresponding
collections over $k$. Then, a partial Benders reformulation with
respect to $\hat E$ and $\hat F$, denoted by $\mathbf{BD-MIBP}(\hat
E,\hat F)$, is a relaxation of $\mathbf{BD-MIBP}$ (as well as
$\mathbf{CC-MIBP}$), and its optimal value yields a  lower
bound.  \hfill $\square$
\end{prop}
It is clear that any feasible solution to $\mathbf{CC-MIBP}$ provides an
upper bound. So, iteratively including constraints from extreme
points and rays, i.e., those in \eqref{eq_xpoints} and
\eqref{eq_xrays} (which are often referred to as Benders optimality
and feasibility cuts respectively), as cutting planes could help us
generate better lower and upper bounds, yielding a variant of
 Benders decomposition method based on our bilinear Benders
 reformulation.

Note that $\mathbf{DP}_k$ is finitely optimal or unbounded for all
$k$. Hence, if CCMP is unbounded, it must be resulted from
$\min\{\mathbf{cx}: \mathbf{x\in X}\}=-\infty$. To avoid such
situation, which may cause Benders decomposition method enumerate all extreme
points and rays before the observation of unboundedness, we assume
that $\min\{\mathbf{cx}: \mathbf{x\in X}\}$ is finite. Let
$\mathbf{CC-MIBP}(\mathbf{x}^0,\mathbf{z}^0)$ represent the
formulation with $\mathbf{x}$ and $\mathbf{z}$ fixed to $
\mathbf{x}^0$ and $\mathbf{z}^0$ respectively, $LB$ and $UB$ be the
current lower and upper bounds, and $e$ be the optimality tolerance.
Our procedure, to which we call (bilinear) Benders decomposition, is
as follows.

\vspace{10pt} \noindent\textbf{\underline{Basic (Bilinear) Benders
Decomposition Method}}
\begin{description}
\item[\underline{Step 1.} -- Initialization:] Set $LB= - \infty$ and $UB=
+\infty$. Set $\hat E=\hat F=\emptyset$.
\item[\underline{Step 2.} -- Iterative Steps:] \textrm{}
\begin{description}
    \item[2.a] Compute (linearized) master problem $\mathbf{BD-MIBP}(\hat E,\hat
    F)$. \\
    $(i)$  If it is
    infeasible, terminate. Otherwise, obtain its optimal value, $V_l$, and
    an optimal solution $(\mathbf{x}^0,\mathbf{z}^0)$. \\
    $(ii)$ Update $LB= V_l$.
    \item [2.b] For all $k$ such that  $z^0_k=0$,\\
    $(i)$ compute  subproblem $\mathbf{DP}_k$ defined in
    (\ref{eq_sub_obj}-\ref{eq_sub_var})  to obtain
    an optimal solution $\mu^0$ if bounded, or derive an  extreme ray
    $\upsilon^0$ if unbounded.\\
    $(ii)$ update $\hat E^k=\hat E^k\cup\{\mu^0\}$ or $\hat F^k=\hat
    F^k\cup\{\upsilon^0\}$, i.e.,  supply
     Benders cuts \eqref{eq_xpoints} or
\eqref{eq_xrays}  to master problem $\mathbf{BD-MIBP}(\hat E,\hat
F)$.
    \item[2.c] Solve $\mathbf{CC-MIBP}(\mathbf{x}^0,\mathbf{z}^0)$ and
     obtain its optimal
    value $V_u$. Update $UB = \min\{UB, V_u\}$.
\end{description}
\item[\underline{Step 3.} -- Stopping Condition:] \textrm{}\\
 If $\frac{UB-LB}{LB}\leq e$, terminate with a solution associated with $UB$.
Otherwise, go to Step 2. \hfill $\square$
\end{description}
The aforemention algorithm is very close to the standard Benders
decomposition method applied to solve SP, except that Benders cuts are
only generated for scenario/subproblem with $z^0_k=0$. We denote
this basic implementation as BD0.
Actually, such a restriction is not necessary. So, as an alternative
approach, we generate Benders cuts from all scenarios/subproblems
regardless of  $z^0_k$'s value. We denote it as BD1. As BD1 has almost the same
simple structure and complexity as BD0, unless otherwise stated, we
call both of them basic Benders decomposition method. Indeed, the
minor change made in BD1 leads to  a significant computational
improvement (see results in Section \ref{Sect_numerical}). Note that
one possible termination scenario of this Benders decomposition method is the
detection of infeasibility specified in Step 2.a.


Next, we present a simple numerical study on a few SP and CCMP
instances, which helps us appreciate the connection between SP and
CCMP, and provides ideas on improving the aforementioned Benders
decomposition procedure in solving CCMP.

\subsection{Observations and Insights on Computing SP and CCMP}
\label{subsect_simplenumerical} We present a rather small-scale
numerical study that involves computing SP formulation,
$\mathbf{CC-bigM}$ formulation, and the linearized
$\mathbf{CC-MIBP}$ formulation. Our objective is to learn the
connection and the difference between SP and CCMP, the impact of the
chance constraint, and to develop a basic understanding on the
performance of Benders decomposition for these two types of
problems. We simply generate one instance using each of six instance
generation combinations described in Section \ref{Sect_numerical}.
In order to highlight the impact of chance constraint and to
minimize the computational challenges from other factors, we
intentionally select instances whose SP formulations can be solved
by professional solver CPLEX in a short time (e.g., less than 60
seconds). Overall, we have 6 instances for binary and general
integer $\mathbf{x}$, respectively. Our computation environment and
related parameters are described in Section
\ref{Sect_numerical}.

Our results are presented in Table
\ref{tab:sample_binary_x}-\ref{tab:sample_general_x}. In those
tables, CPX denotes using CPLEX to compute MIP formulation, BD
denotes Benders decomposition method, CPX-bigM represents computing
$\mathbf{CC-bigM}$ by CPLEX, CPX-MIBP represents computing
linearized $\mathbf{CC-MIBP}$ by CPLEX. If one instance can not be
solved by some method due to time limit, which is set to 3,600 seconds, or memory issue,
we use ``T'' or ``M'' to represent that entry, respectively, and report
available gap information. Row ``avg.: XMS''
presents the average values over instances whose MIP (big-M, if
applicable) formulations can be optimally solved by  CPX.
To have a
fair comparison between SP and CCMP under Benders decomposition
method, we select BD1 to compute CCMP, noting that both BD1 and
Benders decomposition in SP generate Benders cuts from all
scenarios.  In our computational study, big-M is set to $10^5$,
which is also used as variable upper bound (except binary variables)
in McCormick linearization operations.

\begin{table}[htbp]
  \centering
  \caption{Instances with Binary $\mathbf{x}$ }
    \scalebox{0.75}{\begin{tabular}{|c|cc|ccc|cc|cc|ccc|}
    \hline
    \multirow{3}[0]{*}{\#} & \multicolumn{5}{c|}{\textbf{SP}} & \multicolumn{7}{c|}{\textbf{CCMP}} \\
    \cline{2-13}
          &\multicolumn{2}{c|}{CPX}   & \multicolumn{3}{c|}{BD}  & \multicolumn{2}{c|}{CPX-bigM}
             &\multicolumn{2}{c|}{CPX-MIBP} & \multicolumn{3}{c|}{BD}\\
          & sec.   & g(\%) & itr.  & sec.   &g(\%) & sec. &g(\%)  & sec.  &g(\%) & itr.  & sec.  &g(\%)\\
          \hline
    B1    & 6  &   & 3     & 9  &   & 873 &     & 470 &    & 4     & 53 &   \\
    B2    & 8  &   & 4     & 15 &   & T &   21.27   & 969 &    & 5     & 112  &  \\
    B3    & 18 &   & 3     & 26 &   & 1835 &    & 396 &   & 2     & 36 &    \\
    B4    & 7 &    & 4     & 15 &   & 2081 &      & 906 &   & 4     & 58 &  \\
    B5    & 15 &   & 5     & 25 &   & T&   4.65  & 914 &   & 4     & 88  &   \\
    B6    & 33 &   & 5     & 103 &   & T&   11.25   & T &   11.29 & 5     & 438 &   \\
    \hline
    avg.: XMS  & 15 &   & 4.00  & 32 &   & 1596 &     & 591 &  & 3.33  & 49 &   \\
    \hline
    \end{tabular}%
  \label{tab:sample_binary_x}} %
  \centering
  \caption{Instances with General Integer $\mathbf{x}$ }
    \scalebox{0.75}{\begin{tabular}{|c|cc|ccc|cc|cc|ccc|}
    \hline
    \multirow{3}[0]{*}{\#} & \multicolumn{5}{c|}{\textbf{SP}} & \multicolumn{7}{c|}{\textbf{CCMP}} \\
    \cline{2-13}
          &\multicolumn{2}{c|}{CPX}   & \multicolumn{3}{c|}{BD}  & \multicolumn{2}{c|}{CPX-bigM}
          &\multicolumn{2}{c|}{CPX-MIBP} & \multicolumn{3}{c|}{BD}   \\
          & sec.   & g(\%) & itr.  & sec.   &g(\%) & sec. &g(\%)   & sec.  &g(\%) & itr.  & sec.  &g(\%) \\
          \hline
    I1    & 7 &    & 5     & 23 &   & 2050  &      & 611 &   & 3     & 64 &  \\
    I2    & 17 &   & 5     & 25 &   & 257 &      & 348 &   & 2     & 17 &  \\
    I3    & 26 &   & 6     & 106 &   & 1913  &      & 1229 &    & 2     & 72 &   \\
    I4    & 33 &   & 6     & 110 &   & T &55.9     & T &NA   & 6     & T &14.65   \\
    I5    & 59 &   & 4     & 27  &    & T &0.7   & 698 &   & 2  & 35 &     \\
    I6    & 27 &   & 8     & 183 &     & T  &55.16        & T  &13.24   & 5     & T &9.05  \\
    \hline
    avg.: XMS & 28 &   & 5.67  & 79   &     & 1407 &      & 730 &   & 2.33  & 51 &    \\
    \hline
    \end{tabular}}%
  \label{tab:sample_general_x}%
\end{table}%

  Clearly, Table \ref{tab:sample_binary_x}-\ref{tab:sample_general_x} confirm
  the well-known result that chance constrained models are more challenging
  than their SP counterparts, which is more prominent when  solver
  CPLEX is adopted as the solution method. Note that a single chance
  constraint  could easily incur hundreds of times more computational
  expenses. Through  a closer study on Table \ref{tab:sample_binary_x}-\ref{tab:sample_general_x}, we have a few  more interesting observations
  that have not been fully discussed in previous study.

  \begin{enumerate}
  \item From the point of view of the solver,
  the popular $\mathbf{CC-bigM}$ formulation of CCMP is much more difficult to  compute
  than  the  (linearized) bilinear formulation $\mathbf{CC-MIBP}$. Among all
  12 instances,  only 6 of them can be solved if they are
  represented  by $\mathbf{CC-bigM}$ formulation. Nevertheless, if we
  adopt the (linearized) $\mathbf{CC-MIBP}$, 9 out of them can be
  solved to the optimality. Indeed, the bilinear
  formulation typically can be  2 to 4 times faster than the other
  one. Hence, we believe that the bilinear
  formulation could be both theoretically stronger and
   computationally more friendly. It definitely deserves a further
  study.
  \item Benders decomposition method, implemented in the presented bilinear form,
  actually displays a very strong
  capability to compute CCMP, which is different from an understanding made in
  \cite{tanner2010iis} that it is
  not a good method.
  Unlike the    solver, which is severely
  affected by the single chance constraint in every CCMP instance, Benders method often can
  compute optimal solutions with reasonably more time. According to Table
  \ref{tab:sample_binary_x}, Benders method on average produces an optimal solution with
  3-4 times more computational expenses when SP is converted into CCMP. Nevertheless,
   the professional solver
  could fail to derive optimal solutions even with hundreds of times
  more computational expenses.  Hence, it can be seen that  Bender decomposition method
  is much  more robust to the computational challenge caused by the
  chance constraint.
  \item We further notice two more non-trivial points from Table \ref{tab:sample_binary_x}-\ref{tab:sample_general_x} that
  support the effectiveness of Bender
  decomposition method.
  One
  is that the performance benchmark between Benders decomposition method and the
  professional solver demonstrates an opposite behavior in SP and
  CCMP. In SP instances, the solver is typically more
  efficient.
  In CCMP instances, however, Benders decomposition method  generally performs
  an order of magnitude faster than the solver. Another one is that, in terms of iterations,
  CCMP might not be more challenging than SP.
  For some instances, CCMP formulation may just need one
  more iterations, e.g.,
  instances B1 and B2, compared to their SP counterparts.
  Nevertheless, for some other instances, CCMP formulation can be computed with
  noticeably less iterations, e.g., instances I2 and I3. We think that such observation is counterintuitive
  and has not been reported in any previous study. It actually indicates that
  Benders decomposition method, in the presented bilinear form, has a capability to quickly identify their optimal sets of responsive scenarios from
  their scenario pools and derive corresponding solutions. Hence,
  from these two observations, we believe that
  Benders decomposition method is probably more
  appropriate to compute CCMP than to compute SP, given that it can effectively
  deal with the combinatorial structure implied in the chance constraint.

  \end{enumerate}
  Besides using the presented basic
  Benders decomposition method, we can design new enhancement strategies
  or incorporate  existing ones developed for solving  SP to further improve
  computational performance. Next, we present a few improvement
  techniques that make use of the
  (possible) easiness of the underlying SP formulation or CCMP's structural
  properties.

\subsection{Enhancement Strategies of Benders Decomposition}
\label{subsect_enhancement}
Results in Section \ref{subsect_simplenumerical} motivate us to make
use of the rather easier SP model to efficiently compute, exactly or
approximately, CCMP. Following this line, we present a couple of
strategies that can be incorporated into our basic Benders
decomposition method in solving CCMP.


\noindent \underline{$(i)$ SP based initialization}: One straightforward enhancement strategy is to ignore the chance constraint
and simply to solve the corresponding SP model by Benders
decomposition method to obtain a set of Benders cuts. According to
Proposition \ref{prop_validity_Benders}, those Benders cuts, as in \eqref{eq_xpoints} and
\eqref{eq_xrays},
can be employed to initialize $\hat E$ and $\hat F$ in Step 1 of the basic Benders
decomposition method. We refer to this strategy as SP based initialization.

\noindent \underline{$(ii)$ Small-M based initialization}:  Actually, we can directly treat  $\mathbf{CC-bigM}$ formulation
of CCMP as an SP problem, where both $\mathbf{x}$ and $\mathbf{z}$
are first stage variables, and compute it by Benders method.
Nevertheless, the challenge from big-M coefficients of $\mathbf{z}$
remains in computing the resulting master problem. Different from
the popular idea of big-M concept, we can adopt smaller coefficients
for $\mathbf{z}$, to which we denote small-M coefficients, to serve
the initialization purpose. Noting that big-M coefficients are
larger than small-M coefficients, the next result follows naturally.

\begin{prop}
\label{prop_small_M} The popular $\mathbf{CC-bigM}$ formulation of CCMP
is a relaxation of the same formulation with small-M coefficients,
i.e., the \emph{small-M MIP formulation}.  Therefore, any feasible
solution to the small-M MIP formulation is feasible to $\mathbf{CC-bigM}$, and therefore
feasible to CCMP model.
\end{prop}

It can be seen that the small-M MIP formulation is to
approximate CCMP in the primal space, which leads to feasible
solutions and strong upper bounds. We can also derive dual
information by treating it as a regular SP, whose first stage
variables are $\mathbf{x}$ and $\mathbf{z}$, and computing it by
Benders decomposition method. Then, similar to SP based
initialization, the resulting Benders cuts can be employed for
initialization. We refer to this strategy as small-M based
initialization. As the small-M MIP formulation resembles
$\mathbf{CC-bigM}$ more than SP does, those Benders cuts probably
should be more effective than those obtained in SP based initialization.

\noindent \textbf{Remark:} We note that small-M MIP formulation might
provide a flexible framework to analyze and compute CCMP. When
$M=0$, small-M MIP formulation simply reduces to SP model, which defines
a core set inside $\mathbf{X}$ that is feasible to all scenarios. It
is worth mentioning that core set is convex. When $M$
increases, small-M MIP formulation approaches CCMP, which also expands
that core set. Finally, when small-M is sufficiently large, it
becomes the popular $\mathbf{CC-bigM}$ formulation for CCMP and that set will become the
feasible set of CCMP, which is typically non-convex. Actually, this
discussion of small-M MIP formulation might reveal the transition
between SP and CCMP, noting that SP is easy to solve and
$\mathbf{CC-bigM}$ is highly difficult to compute. Hence, a trade
off between the computational expense and the solution quality could
be achieved by adjusting the value of small-M coefficients.
Furthermore, we can interpret that both SP and small-M MIP
formulations are to compute (approximately) an optimal
solution of CCMP and the presented bilinear Benders decomposition is
to perform verification of the optimality of that solution.

\noindent \underline{$(iii)$  Benders decomposition with integer cuts}:   Consider a given
$\mathbf{z}^0$. Let $\mathbf{CC-MIBP}(\cdot,\mathbf{z}^0)$
represent CCMP formulation with $\mathbf{z}=\mathbf{z}^0$, and
denote its optimal value by $\tilde V_u(\mathbf{z}^0)$. Note that
$\mathbf{CC-MIBP}(\cdot,\mathbf{z}^0)$ is essentially an SP problem
that could be easy to compute. Indeed, if $\mathbf{z}^0$ has been
evaluated, we can remove it from  master problem
to reduce the computational time. It can be achieved by adding the
following \textit{integer cut}
$$|\mathbf{z}-\mathbf{z}^0|\geq 1$$ into master problem. Given the binary property of
$\mathbf{z}$, we can easily simplify it. Let set $\mathbf{K}_1$ be
collection of indices of $z^0_k=1$ and $\mathbf{K}_0$ be the
complement set. The integer cut can be formulated as
\begin{eqnarray}
\label{eq_integr_cut}\sum_{k\in \mathbf{K}_1}z_k - \sum_{k\in
\mathbf{K}_0}z_k\leq |\mathbf{K}_1|-1.
\end{eqnarray}

Accordingly, we modify the basic Benders decomposition procedure to
generate stronger bounds. Specifically, let $\mathbf{BD-MIBP}(\hat
E,\hat F| \mathbb{C})$ represent master problem
$\mathbf{BD-MIBP}(\hat E,\hat F)$ subject to additional conditions
$\mathbb{C}$. In Step 1, we include the initialization
$\mathbb{C}=\emptyset$. In Step 2.a, we replace
$\mathbf{BD-MIBP}(\hat E,\hat F)$ by $\mathbf{BD-MIBP}(\hat E,\hat
F| \mathbb{C})$, and modify Step 2.a.($i$) as ``If it is
    infeasible, terminate (with a solution associated with $UB$ if
    it exists). $\dots$'' In Step 2.c, we replace
$\mathbf{CC-MIBP}(\mathbf{x}^0,\mathbf{z}^0)$
 by $\mathbf{CC-MIBP}(\cdot,\mathbf{z}^0)$.
 It is clear that those changes yield new bounds that are stronger than those
 produced in the basic Benders decomposition procedure. We refer to this
improvement as Benders decomposition with integer cuts.

\noindent \underline{$(iv)$  Benders decomposition with Jensen's inequality}: For one type of instances where $\mathbf{H}_k=\mathbf{H}$ and
 $\mathbf{f}_k=\mathbf{f}$ for all $k$, 
 Jensen's inequality derived in Section \ref{sect_Jensen} can be used to improve
 the basic Benders decomposition method. 
 The
next result directly follows Theorem \ref{thm_Jensens} and Corollary \ref{cor_jensens}.
 \begin{prop}
 \label{prop_jensen_computation} Let $E\eta$ be a variable representing the expected recourse cost from responsive scenarios.
The following inequalities
 are
 valid for $\mathbf{CC-MIBP}$:
 \begin{eqnarray}
  \label{eq_Jensen_BD_original}
E\eta &\geq& (1-\sum_{k=1}^K\pi_kz_k) \mathbf{f}\mathbf{\overline
 y}, \\
\mathbf{H}\mathbf{\overline y} &\geq& \frac{\sum_{k=1}^K\pi_k(1-z_k)}{ (1-\sum_{k=1}^K\pi_kz_k)}
 (\mathbf{h}_k-\mathbf{G}_k\mathbf{x}). \label{eq_Jensen_BD_original2}
 \end{eqnarray}
 where  $\mathbf{\overline y} \in \mathbb{R}^m_+$ are newly
 introduced variables. When $\pi_k=\frac 1 K$ for all $k$ and $\mathbf{f}\geq \mathbf{0}$, they can be simplified as
  \begin{eqnarray}
  \label{eq_Jensen_equalPI}
 E\eta  &\geq& \frac{K-L}{K}\mathbf{f}\mathbf{\overline
 y}, \\
\mathbf{H}\mathbf{\overline y} &\geq& \sum_{k=1}^K\frac{1-z_k}{K-L} 
 (\mathbf{h}_k-\mathbf{G}_k\mathbf{x}).     \label{eq_Jensen_equalPI2}
 \end{eqnarray}
 In particular, if $\mathbf{G}_k=\mathbf{G}$ for all $k$,  \eqref{eq_Jensen_equalPI2} reduces to the following linear inequality 
 \begin{eqnarray}
 \label{eq_Jensen_linear}
 \mathbf{H}\mathbf{\overline y} \geq \sum_{k=1}^K\frac{1-z_k}{K-L}\mathbf{h}_k- \mathbf{Gx}.
 \end{eqnarray} \hfill $\square$
 \end{prop}
 Within Benders decomposition scheme, the aforementioned variables and inequalities
 can be included in master problem $\mathbf{BD-MIBP}(\hat E,\hat
 F)$ from the beginning. Again, given the binary property of $z_k$,
 bilinear terms, if exist, 
 can be linearized using McCormick linearization method.  The augmented master problem with
 those variables and constraints will have a stronger lower bound
 that may lead to faster convergence.  As demonstrated in a few
 SP applications 
\cite{erdogan2013dynamic,batun2011operating}, Benders Decomposition
method strengthened by such strategy significantly outperforms other
Benders Decomposition variants, including an implementation within a
Branch-and-Cut framework.

 We observe that including new variables and constraints, especially constraints in \eqref{eq_Jensen_BD_original2}, \eqref{eq_Jensen_equalPI2} 
 or \eqref{eq_Jensen_linear} have a clear impact on the complexity of master problem. So, instead of applying
 the aforementioned Jensen's inequality, we can 
 adopt some relaxation of Jensen's inequality within master problem. Next, we present one such relaxation 
 for the case where $\mathbf{G}_k=\mathbf{G}$ and $\pi_k=\frac 1K$ for all $k$, and $\mathbf{f}\geq \mathbf{0}$.
 
Let  $h_{k,i}$  be the $i^\textrm{th}$ component of $\mathbf{h}_k$.  For a fixed $i$, we  sort $h_{k,i}$ from small to large  and  use $k_s$ to return the original (scenario) index of the $k^\textrm{th}$ one in this sorted sequence. Accordingly, we compute the conditional mean $\overline h_i$ based on the smallest $(K-L)$ ones of this sequnce, i.e., $\overline {h}_i = \frac{\sum_{k=1}^{K-L}h_{k_s,i}}{K-L}$.  We then let $\mathbf{ \overline h}$ be the vector of $\overline {h}_i$ over all $i$,  which leads to the next result. 

\begin{cor}
\label{cor_Jensen_relax_linear}
The inequality \ 
\begin{eqnarray}
\label{eq_jensen_relaxation}
\mathbf{H}\mathbf{\overline y} \geq \mathbf{\overline h} -\mathbf{G}\mathbf{x}
\end{eqnarray} is a relaxation to  \eqref{eq_Jensen_linear}.  So, it can be used to replace \eqref{eq_Jensen_linear}  to obtain a valid relaxed Jensen's inequality for $\mathbf{CC-MIBP}$. \hfill $\square$
\end{cor}
\noindent Because $\varepsilon$ is typically small, which suggests $(K-L)$ is close to $K$,   \eqref{eq_jensen_relaxation}  could provide a strong lower bound to support Benders decomposition.  Moreover, noting that this relaxed one is of a simple structure same as that of the classical one in     \eqref{eq_jensen_SP}  for SP,  we believe that it could be computationally more   friendly than the exact one in \eqref{eq_Jensen_linear}.

We mention that all the presented enhancement strategies are
designed  explicitly based on CCMP's structure or its Jensen's inequality. Indeed, the bilinear  
 Benders decomposition method presented in this section is
 rather a general framework,
 which can easily be modified to incorporate almost all types of existing Benders
 enhancement strategies, such as
 Pareto optimal cut
 \cite{magnanti1981accelerating,papadakos2008practical},
 multicut aggregation
 \cite{birge1988multicut,trukhanov2010adaptive}, and maximum feasible subsystem cut generation
 \cite{saharidis2010improving}. Given the enormous amount of research
results on Benders decomposition and SP, it is definitely worth
investigating  the integration or customization of those results
 within the presented bilinear Benders variant to deal with
 complicated real problems. In Section \ref{Sect_numerical}, a preliminary
 study  on random
 CCMP instances and instances of chance constrained operating room scheduling problem is presented to appreciate the basic Benders decomposition method and
 enhancement strategies described in this section.

\section{Computational Experiments}\label{Sect_numerical}

\subsection{Computing Environment and Solution Methods} 
Our computation is
made through Concert Technology of CPLEX 12.4, a state-of-the-art
professional MIP solver on a Dell Optiplex 760 desktop computer (Intel Core 2 Duo
CPU, 3.0GHz, 3.25GB of RAM) with Windows XP platform. Benders
decomposition algorithms are implemented in C++ and CPLEX in the
same environment. Table \ref{alg_num} summarizes our twelve different computing methods, including eleven types of Benders decomposition implementations. Among them, we mention that CPX represents using
CPLEX to compute the currently popular $\mathbf{CC-bigM}$
formulation with big-M set to $10^5$. BD2 denotes a BD1
implementation that generates Pareto optimal Benders cuts
\cite{magnanti1981accelerating,papadakos2008practical}, whose
effectiveness has been observed in many applications to address SP
problems. For BD5-BD8, small-M based initialization is performed
with small-M set to 1000. In BD7 and BD8, we adopt a \emph{strongest
cuts only}  strategy,  which just employs the strongest Benders cut
with respect to the current $\mathbf{x}^0$ in each scenario to
perform initialization and ignores other Benders cuts derived from
computing small-M MIP formulation. In BD8, to generate a pool of strong
Benders cut for initialization, we further implement Pareto optimal
Benders cuts when computing small-M MIP formulation. We mention that, on top of each of  those  9 different
Benders decomposition implementations,  Jensen's inequality or the relaxed Jensen's inequality described in Section \ref{subsect_enhancement} can be integrated into its master problem.  To evaluate the benefit of adding those inequalities, we select BD1 as the basis and study the integration with Jensen's inequality or  the relaxed one, which are denoted by BD1J and BD1RJ respectively.

In our computational study, the time limit is set to 3,600 seconds.
Also, when computing $\mathbf{CC-bigM}$ formulation and master
problem of Benders decomposition, the optimality gap $e$ is set to
$0.005$. For subproblem $\mathbf{DP}_k$, the optimality gap $e$ is
set to $10^{-4}$. If initialization operations are involved, we set
the optimality gap of master problem in the initialization stage to
0.02 and restrict the total time for the initialization stage less
than 500 seconds.

\begin{table}[H]
  \centering
  \caption{Description of CPLEX and Benders Decomposition Methods}
    \scalebox{0.8}{\begin{tabular}{|c|c|}
    \addlinespace
    \hline
     Method & Description \\
    \hline
    CPX     &  \multicolumn{1}{l|}{CPLEX with $\mathbf{CC-bigM}$ formulation} \\
    BD0     &  \multicolumn{1}{l|}{Basic method generating cuts from scenarios with $z_k=0$ only} \\

    BD1     &  \multicolumn{1}{l|}{Basic method generating cuts from all scenarios} \\

    BD2     &  \multicolumn{1}{l|}{BD1 with Pareto optimal cuts} \\

    BD3    & \multicolumn{1}{l|}{BD1 with SP based initialization} \\

     BD4 &  \multicolumn{1}{l|}{BD1 with SP based initialization and integer cuts}\\

     BD5     & \multicolumn{1}{l|}{BD1 with small-M based initialization} \\

     BD6    &  \multicolumn{1}{l|}{BD1 with small-M based initialization and integer cuts} \\

     BD7   & \multicolumn{1}{l|}{BD1 with  small-M based initialization (\emph{strongest cuts only}) and integer cuts} \\

     BD8     & \multicolumn{1}{l|}{BD1 with  small-M based initialization (\emph{strongest cuts only}), Pareto} \\
                                      & \multicolumn{1}{l|}{optimal cuts in computing small-M MIP formulation, and integer cuts}  \\
     BD1J     & \multicolumn{1}{l|}{BD1 with Jensen's inequality listed in Proposition \ref{prop_jensen_computation}} \\
     BD1RJ     & \multicolumn{1}{l|}{BD1 with  relaxed Jensen's inequality listed in Corollary  \ref{cor_Jensen_relax_linear} } \\                                      
    \hline
\end{tabular}}
  \label{alg_num}%
\end{table}%

\subsection{Test Beds}
Our numerical study primarily involves unstructured
instances that are generated randomly. Specifically, let $K$, $I_1$, $I_2$
$n$, and $m$ denote the number of scenarios, the number of
constraints in (\ref{mibp_eq_first}), the number of constraints in
(\ref{mibp_scen_const}), the dimension of $\mathbf{x}$, and the
dimension of $\mathbf{y}_k$, respectively. In our study, we consider
two setups of ($I_1$, $I_2$, $n$, $m$), i.e., $T_1=(10,30,20,40)$ and
$T_2=(20,50,30,70)$, and three different $K$ values, i.e., $250,
500$ and $1000$.  For each of six resulting combinations, we
generate five instances where their coefficients are randomly
selected from the ranges described in Table
\ref{tab:parameter_range}.
\begin{table}[H]
  \centering
  \caption{Data Ranges}
    \scalebox{0.8}{\begin{tabular}{|c|l|l|}
    \addlinespace
    \hline
    \textrm{Parameter} & \textrm{Dimension} & \textrm{Range} \\
    \hline
    $A$     & $I_1 \times n$   & [-25,25] \\
    $b$     & $I_1$     & [-50,50] \\
    $c$     & $n$     & [100,300] \\
    $G_k$  & $I_2 \times n$   & first 40\%$I_2$ rows: [0,10], other rows: 0 \\
    $H_k$  & $I_2 \times m$   & first 40\%$I_2$ rows: [-3, 0], other rows: [0,3] \\
    $h_k$  & $I_2$     & first 40\%$I_2$ elements: [-35, 0], other elements: [-25,100] \\
    $f_k$     & $m$     & [5,10] \\
    \hline
    \end{tabular}}%
  \label{tab:parameter_range}%
\end{table}%

Following those specifications, we generate one set of instances with
$\mathbf{x}$ being binary variables and another set with
$\mathbf{x}$ being general integer variables, to study algorithm
performance with different types of $\mathbf{x}$. For $\mathbf{x}$
being general integer variables, we set their upper bounds to 500,
which is also used in linearization operations in Benders
decomposition.    Overall, we have $30\times 2$
instances. Note that it is not necessary to require $\mathbf{x}$
to be integer variables in our algorithm implementation. 
In order to
gain insights on  benchmarking Benders decomposition method and the
professional solver CPLEX, we primarily consider random instances of
a moderate difficulty (i.e., the computational time by CPLEX is more
than 500 seconds and has less than 15\% optimality gap before time
limit when under $\varepsilon=0.1$).  To evaluate the solution capability on more challenging instances, 
we also generate 30 difficult ones whose first stage $\mathbf{x}$ are general
integer and  have a gap  $\geq$15\% by CPLEX before time limit. 
  
 It has been demonstrated in a couple of stochastic planning or scheduling problems \cite{erdogan2013dynamic,batun2011operating} that Jensen's inequality could be very beneficial to Benders decomposition when its basic version is inefficient. In order to appreciate its advantage in solving 
chance constrained programs,  we extend the stochastic operating room (OR) scheduling model presented in \cite{batun2011operating} to build its chance constrained formulation. Appendix of this paper presents this  formulation, along with the description of associated parameters and decision variables.  Note that the structure of this formulation enables us to incorporate Jensen's inequality presented in Proposition \ref{prop_jensen_computation} or  the relaxed one presented in Corollary \ref{cor_Jensen_relax_linear} into master problem. Then, random instances  of this particular application are generated according to the specifications made in \cite{batun2011operating}:  8 surgeries with index $i$ or $j$, 2 surgeons with index $k$, 2 operating rooms with index $q$ or $r$, and 9 hours per day ($L$ = 540).  Surgery list of surgeon 1  is \{1,2,3,4\} and that of surgeon 2 is \{5,6,7,8\}. The fixed cost $c^f$, overtime cost $c^o$, and idle time cost $c^S$  are 4437, 12.37, and 17.748, respectively. Surgeon turnover time $s^S$ and operating room turnover time $s^R$ are 0 and 30, respectively. The big $M$ value is set to 2500.  Durations of different stages, i.e., preincision, incision and postincision,  of an operation are random.  In our numerical study, we consider two situations with different random levels. In the first situation, we employ $K$ = 100 scenarios with equal probabilities to capture the uncertainty, where $pre_i(\omega)$ and $post_i(\omega)$  are randomly selected from $[32, 56]$ and $p_i(\omega)$ are randomly selected from  $[50, 150]$, $\forall i$ of a particular scenario. In the second one,  those parameters of each scenario are randomly selected from $[26, 38]$, $[26, 38]$ and $[76, 123]$, respectively. Clearly, the first situtation demonstrates more severe randomness than the second one.  
 Then,  we generate 5 instances for each situation to support our computation.

\subsection{Computational Results of Benders Decomposition Variants} 
\subsubsection{Results of Unstructured Random Instances}
We first present numerical results of  unstructured random instances in 
Table \ref{base_bin} and
\ref{base_int}, with $\varepsilon$ set to 0.1. In those tables, 
column ``obj.'' presents the objective function value of the best
feasible solution. Column ``sec.'' records the computational time in
seconds. If the algorithm is terminated because of time limit, we
label it by ``T''. If it is terminated due to insufficient memory, we
label it by ``M''. Column ``g($\%$)'' keeps the relative gap in
percentage when the algorithm terminates before reaching optimality.
Column ``itr.'' presents the total number of iterations in Benders
decomposition. If any of initialization techniques is applied, the
total number of iterations will be displayed as the number of
iterations in the initialization stage $+$ that number in the
regular computation stage. At the bottom of those tables, we have a
few rows to show the average performances of different computing
methods. Row ``\# solved (S)'' provides the number of instances
solved to optimality. Similarly, row ``\# unsolved (U)'' provides
the number of instances unsolved. Row ``avg. sec.: S'' presents the
average computational time in seconds over instances solved to
optimality. Row ``avg. gap: U'' presents the average relative gap in
percentage over  unsolved instances. As there are a few instances
that do not have gap information, we do not include them when
computing this row. Row ``CPX/BD: XBS'' displays the average of
computational time ratios between CPLEX and BD, over instances solved
by both CPLEX and Benders decomposition.  Based on those two tables, we make a few
observations: 

\noindent$(i)$ Benders decomposition method, even in its basic forms
BD0 and BD1, shows a significantly better performance than that of
the professional solver CPLEX. With the most effective
enhancement strategies, Benders decomposition optimally solves 95\% of all
instances while CPLEX fails to derive optimal solutions for 50\% of
those instances. For those that can be computed by both CPLEX and
Benders decomposition, according to ``CPX/BD: XBS'',  the latter one
generally is faster by an order of magnitude. Indeed, given that our
testing instances are unstructured, we can anticipate that such
improvement can be more prominent in well-structured real
applications  where Benders decomposition can be highly customized. \\
$(ii)$ Although BD1 is obtained by making a simple modification on
BD0, its performance is clear superior to BD0. Enhancement
strategies, especially SP based and small-M based initialization
techniques, could further improve algorithm performance at a
significant level. For example, when $\mathbf{x}$ are binary, the
best performed method is BD5, which is with small-M based
initialization. Comparing BD1 and BD5, we note that more than 25\%
computational time is reduced. And on those solved by CPLEX, it
perform 22 times faster. When $\mathbf{x}$ are general integers, it
is interesting to note that although BD5 probably works faster on
those solved by CPLEX, BD3, which is with SP based initialization,
performs more stably by solving more instances than BD5.

Other enhancement strategies generally are not as effective as BD3
and BD5. One explanation is that those instances can be computed
with a small number of Benders iterations, which may not provide a
suitable platform to demonstrate the effectiveness of other strategies. Nevertheless,
there may exist some applications where other enhancement
strategies could be very effective.\\
%
$(iii)$ Results on all Benders decomposition implementations with SP or small-M based initialization
empirically confirm our discussion following Proposition
\ref{prop_small_M}. It is often the case that one single
bilinear Benders iteration is needed after the initialization stage
is completed (with a high quality feasible solution). Clearly, that single iteration just serves the
verification purpose to the optimality of that solution.  It is more
obvious in BD5, where small-M based initialization is adopted, than in 
BD3, in both Tables. Such observation indicates that this approximation-verification scheme is of a practical interest where
the optimality tolerance is an input parameter from system operator. By adjusting the value of $M$ and using
this approximation-verification scheme, we will be able to produce a solution satisfying the optimality requirement with
a tolerable computational expense for large-scale practical problems.

  To study algorithm performance under a different risk tolerance
   level, we recompute all those instances with $\varepsilon=0.05$. BD5 and BD3 are
   selected as computing methods for instances with binary $\mathbf{x}$ and with general
   integer $\mathbf{x}$, respectively. Results are presented in Table \ref{small_e}.  There are two points worth our attention.
    First, under $\varepsilon=0.05$, those instances become significantly easier for CPX.
    Consequently, the benefit of using Benders
    decomposition over CPLEX is less attractive. One explanation is
    that with $\varepsilon=0.05$,  scenario selection is of less
    freedom, which suggests identifying optimal responsive scenarios is not
    challenging. So, those instances become more suitable for CPLEX. 
     Another point is the strong capability of detecting infeasibility
     demonstrated by  Benders decomposition method. In Table \ref{small_e},
     there are two instances with $K=250$ and $T_2$ that are
     labeled ``inf'', i.e., they are infeasible. On the one hand, it takes
     non-trivial computational expense for CPLEX to report infeasibility. On the other
     hand, BD5 identifies their infeasibility by a couple of
     iterations with much less computational time.  As we mention in
     Section \ref{sect_formulation_assum}, such strong capability is
     of a critical value for practitioners to perform sensitivity analysis with respect to
     $\varepsilon$.

    We further explore the computational capability of Benders decomposition method on 30 difficult
    instances. Table
    \ref{hard}  reports results where BD1, BD3, BD5 are adopted as computing methods. Again, all those Benders decomposition variants demonstrate much more powerful
    computing capabilities than CPLEX. Among instances that can be computed to optimality, BD5 computes much faster than BD1 and BD3 while BD3 performs
    more stably by solving more  instances. For those whose gap information are available before time limit, 
    those Benders decomposition variants always ensure a solution that has a significantly smaller
    gap than that produced by CPLEX. Among those Benders decomposition variants, BD5 is most effective on reducing gaps.  Overall, from the computation time and gap information of those instances, we generally believe that small-M based initialization probably is a more effective strategy than SP based initialization. 
    
\subsubsection{Results of Chance Constrained OR Scheduling Problem}    
    For instances of chance  constrained OR scheduling problem,  we test one basic Benders decomposition variant,  i.e., BD1, and its integration with Jensen's inequalities, i.e., BD1J and BD1RJ. Table \ref{JI} presents their computational performances with $\varepsilon$ equal to 0, 0.05 or 0.1.  Noting that when $\varepsilon = 0$,  chance constrained model reduces to its SP counterpart, BD1 reduces to the classical Benders decomposition method in computing SP, and BD1J and BD1JR are identical to BD1.  To gain insights of BD1J and BD1RJ in computing chance constrained instances,  the last row of Table \ref{JI} only summarize  results of those with  $\varepsilon$ equal to 0.05 or 0.1.   Same as our observations made in Section \ref{subsect_simplenumerical},  although their SP counterparts are rather easy to compute by CPLEX, chance constrained instances are very challenging to solve.  Although BD1 is not an effective solution method 
 for either   SP or chance constrained instances, the benefit from including Jensen's inequalities is drastic. Indeed, such benefit is more obvious in computing chance constrained instances where  BD1J and BD1RJ demonstrate significantly stronger solution capabilities over CPLEX.  Noting that CPLEX typically cannot produce a solution of a small gap for any chance constrained instance,  BD1J or BD1RJ typically can derive solutions with practically acceptable gaps.  Such gap reduction could be in one order of magnitude.  
 
 It actually is interesting to note that BD1J and BD1RJ display different behaviors.  In group 1, which has more severe uncertainties, BD1J generally produces solutions with less gaps. In group 2, which has less randomness,  BD1RJ performs better.  On the one  hand, the results in group 1 can be explained by the fact that when uncertainty gets more severe,  the relaxed Jensen's inequality becomes weaker.  Although it is more computationally friendly with more iterations done within the time limit,  its support to improving the lower bound is not as effective as Jensen's inequality. On the other hand,  in group 2 where less randomness is involved in,  strength of the relaxed Jensen's inequality is close to Jensen's inequality. Given that master problem with Jensen's inequality is often terminated due to its computational complexity,  the advantage of the simple structure of the relaxed one, along with its strength, renders BD1RJ demonstrate  better computational performance over BD1J.  Such different behaviors indicate that   we  need to adopt the appropriate Jense's inequality to achieve the best performance.

\begin{landscape}
\begin{table}[htbp]
  \centering
  \caption{Comparison of Algorithm Performance under $\varepsilon = 0.1$ on Instances with Binary $\mathbf{x}$}
    \scalebox{0.72}{\setlength{\tabcolsep}{0.01em} \begin{tabular}{|c|c|c|lcr|lcr|lcr|lcr|lcr|lcr|lcr|lcr|lcr|lcr|}
    \addlinespace
     \hline
    \multirow{2}[0]{*}{$K$} & \multirow{2}[0]{*}{} & \multirow{2}[0]{*}{} & \multicolumn{3}{|c|}{CPX}  & \multicolumn{3}{c|}{BD0} & \multicolumn{3}{c|}{BD1} & \multicolumn{3}{c|}{BD2} & \multicolumn{3}{c|}{BD3} & \multicolumn{3}{c|}{BD4} & \multicolumn{3}{c|}{BD5}  & \multicolumn{3}{c|}{BD6} & \multicolumn{3}{c|}{BD7}  & \multicolumn{3}{c|}{BD8} \\
  \cline{4-33}
          &       &   & obj.    & sec.   & g(\%) & itr. & sec.  & g(\%) & itr. & sec.  & g(\%) & itr. & sec.  & g(\%) & itr. & sec.  & g(\%) & itr. & sec.  & g(\%) & itr. & sec.  & g(\%) & itr. & sec.  & g(\%) & itr. & sec.  & g(\%) & itr. & sec.  & g(\%) \\
    \hline
    \multirow{10}[0]{*}{250} & \multirow{5}[0]{*}{$T_1$} & 1 & 1893.3   & 873    &    & 4     & 56    &     & 4     & 53    &       & 4     & 92    &    & 3+3   & 63    &       & 3+3   & 87    &   & 4+1   & 35    &         & 4+1   & 45    &       & 4+1   & 43    &            & 4+1   & 35    &  \\
          &       & 2  & 2038.4  & 1684   &   & 3     & 33    &     & 3     & 30    &       & 3     & 46    &    & 3+2   & 35    &        & 3+2   & 29    &    & 3+1   & 22    &         & 3+1   & 38    &       & 3+1   & 34    &          & 3+1   & 54    &  \\
          &       & 3 & 1432.1  & 2081   &    & 4     & 66    &     & 4     & 58    &       & 4     & 103   &   & 4+2   & 65    &          & 4+2   & 80    &    & 4+1   & 35    &        & 4+1   & 47    &       & 4+1   & 43    &          & 4+1   & 36    &  \\
          &       & 4  & 1744.2  & T      & 4.65  & 5     & 115   &   & 4     & 88    &       & 4     & 129   &   & 5+2   & 97    &          & 5+2   & 125   &   & 4+1   & 58    &         & 4+1   & 71    &       & 4+1   & 60    &          & 4+1   & 51    &  \\
          &       & 5  & 1244.9  & 541    &   & 3     & 34    &      & 3     & 32    &       & 3     & 62    &   & 5+2   & 74    &        & 5+2   & 98    &   & 3+1   & 23    &          & 3+1   & 36    &       & 3+1   & 33    &           & 3+1   & 25    &  \\
    \cline{2-3}
          & \multirow{5}[0]{*}{$T_2$} & 1  & 1408.1   & T      & 0.88   & 3     & 69    &   & 2     & 45    &       & 2     & 66    &   & 5+1   & 98    &        & 5+1   & 100   &   & 3+1   & 40    &         & 3+1   & 89    &       & 3+1   & 93    &            & 3+1   & 31    &  \\
          &       & 2   & 3060.8  & T      & 12.3 & 5     & 3252  &   & 3     & 1296  &       & 4     & 2263  &    & 5+3   & 1434  &        & 5+2   & 604   &  & 4+2   & 1229  &           & 4+1   & 381   &       & 4+1   & 408   &          & 4+1   & 467   &  \\
          &       & 3   & 2764.4  & T      & 1.14 & 3     & 189   &   & 2     & 95    &       & 2     & 166   &    & 2+1   & 83    &       & 2+1   & 83    &   & 3+1   & 106   &          & 3+1   & 362   &       & 3+1   & 352   &           & 3+1   & 399   &  \\
          &       & 4   & 3241.7  & T      & 0.72 & 4     & 348   &   & 3     & 259   &       & 3     & 275   &   & 2+2   & 228   &         & 2+2   & 286   &    & 4+1   & 128   &       & 4+1   & 228   &       & 4+1   & 228   &            & 4+1   & 114   &  \\
          &       & 5   & 1768.6  & T      & 7.75 & 5     & 2186  &   & 4     & 247   &       & 4     & 366   &   & 4+2   & 206   &        & 4+2   & 304   &   & 4+1   & 154   &           & 4+1   & 238   &       & 4+1   & 234   &          & 4+1   & 134   &  \\
    \cline{1-3}
    \multirow{10}[0]{*}{500} & \multirow{5}[0]{*}{$T_1$} & 1  & 1361.2   & T      & 11.25 & 5     & 363   &  & 5     & 438   &       & 5     & 573   &   & 5+3   & 496   &         & 5+3   & 609   &   & 5+1   & 196   &         & 5+1   & 233   &       & 5+1   & 216   &           & 5+1   & 208   &  \\
          &       & 2   & 984.3   & 1835  &    & 3     & 62    &     & 2     & 36    &       & 2     & 89    &    & 3+1   & 62    &         & 3+1   & 62    &   & 3+1   & 54    &           & 3+1   & 80    &       & 3+1   & 80    &        & 3+1   & 55    &  \\
          &       & 3   & 809.1  & T      & 0.56  & 3     & 67    &  & 2     & 37    &       & 2     & 56    &    & 3+1   & 58    &        & 3+1   & 57    &    & 3+1   & 54    &      & 3+1   & 81    &       & 3+1   & 92    &             & 3+1   & 168   &  \\
          &       & 4   & 1904.6  & T      & 11.13 & 5     & 374   &  & 5     & 366   &       & 5     & 555   &    & 2+4   & 362   &       & 2+4   & 450   &   & 5+1   & 185   &          & 5+1   & 214   &       & 5+1   & 206   &           & 5+1   & 184   &  \\
          &       & 5   & 390.1  & T      & 0.99  & 3     & 63    &   & 2     & 37    &       & 2     & 91    &    & 3+1   & 55    &        & 3+1   & 55    &   & 3+1   & 55    &         & 3+1   & 74    &       & 3+1   & 74    &           & 3+1   & 54    &  \\
     \cline{2-3}
          & \multirow{5}[0]{*}{$T_2$} & 1  & 1348.4   & T      & 0.88 & 3     & 218   &    & 2     & 168   &       & 2     & 210   &   & 5+1   & 301   &         & 5+1   & 304   &   & 3+1   & 113   &           & 3+1   & 362   &       & 3+1   & 384   &         & 3+1   & 95    &  \\
          &       & 2  & 2282   & 528    &    & 3     & 436   &      & 2     & 135   &       & 2     & 225   &   & 2+1   & 103   &        & 2+1   & 104   &   & 3+1   & 82    &          & 3+1   & 303   &       & 3+1   & 298   &           & 3+1   & 90    &  \\
          &       & 3  & 1736.8   & 3587   &    & 3     & 370   &     & 2     & 107   &       & 2     & 178   &    & 3+1   & 139   &        & 3+1   & 137   &  & 3+1   & 107   &         & 3+1   & 411   &       & 3+1   & 408   &            & 3+1   & 103   &  \\
          &       & 4  & 1301   & T      & 1.54   & 3     & 119   &  & 2     & 88    &       & 2     & 177   &   & 3+1   & 124   &         & 3+1   & 123   &  & 3+1   & 70    &         & 3+1   & 213   &       & 3+1   & 214   &            & 3+1   & 63    &  \\
          &       & 5  & 2448.5   & 1375   &    & 3     & 207   &      & 2     & 286   &       & 2     & 271   &   & 2+1   & 89    &         & 2+1   & 91    &   & 3+1   & 229   &         & 3+1   & 431   &       & 3+1   & 434   &           & 3+1   & 220   &  \\
    \cline{1-3}
    \multirow{10}[0]{*}{1000} & \multirow{5}[0]{*}{$T_1$} & 1   & 1673   & 1256  &   & 3     & 389   &       & 3     & 385   &       & 3     & 739   &   & 3+2   & 439   &         & 3+2   & 854   &   & 3+1   & 242   &          & 3+1   & 380   &       & 3+1   & 360   &          & 3+1   & 243   &  \\
          &       & 2  & 1411.1   & 1610   &    & 3     & 231   &     & 2     & 133   &       & 2     & 335   &   & 3+1   & 308   &         & 3+1   & 308   &    & 3+1   & 191   &          & 3+1   & 298   &       & 3+1   & 298   &         & 3+1   & 181   &  \\
          &       & 3  & 1208.1   & 1249   &    & 3     & 202   &     & 2     & 128   &       & 2     & 352   &   & 3+1   & 203   &         & 3+1   & 188   &    & 3+1   & 185   &       & 3+1   & 266   &       & 3+1   & 270   &            & 3+1   & 181   &  \\
          &       & 4   & 1467.6   & 1644  &   & 3     & 205   &     & 2     & 127   &       & 2     & 201   &    & 3+1   & 297   &        & 3+1   & 290   &   & 3+1   & 175   &           & 3+1   & 271   &       & 3+1   & 268   &         & 3+1   & 607   &  \\
          &       & 5  & 1087.8   & T      & 11.97 & 3     & 249   &    & 3     & 354   &       & 3     & 683   &   & 4+2   & 504   &        & 4+2   & 725   &   & 3+1   & 211   &            & 3+1   & 335   &       & 3+1   & 310   &         & 3+1   & 213   &  \\
    \cline{2-3}
          & \multirow{5}[0]{*}{$T_2$} & 1  & 1436.3   & 2862   &   & 3     & 464   &      & 2     & 347   &       & 2     & 697   &   & 2+1   & 331   &        & 2+1   & 330   &   & 3+1   & 236   &           & 3+1   & 961   &       & 3+1   & 985   &          & 3+1   & 236   &  \\
          &       & 2  & 2136.3   & T      & 7.93  &3 &725 & & 2     & 538   &       & 2     & 907   &   & 2+1   & 413   &        & 2+1   & 425   &    & 3+1   & 307   &         & 3+1   & 1175  &       & 3+1   & 1196  &           & 3+1   & 333   &  \\
          &       & 3  & 1600.6   & 3073   &  & 5     & T     & 37.24    & 5     & T     & 35.59 & 4     & T     & 33.29 & 4+1   & M & 38.01  & 4+1   & M & 38.01 & 4+3   & T     & 37.03 & 4+1         & M & 37    & 4+1   & M & NA     & 4+2   & M & 36.69 \\
          &       & 4   & 2565.7   & 3480  &    & 3     & 994   &      & 2     & 881   &       & 2     & 1005  &   & 2+1   & 831   &        & 2+1   & 818   &   & 3+1   & 366   &           & 3+1   & 1443  &       & 3+1   & 1548  &          & 3+1   & 394   &  \\
          &       & 5    & 2320.7  & 2746  &    & 3     & 373   &     & 2     & 352   &       & 2     & 609   &    & 2+1   & 242   &        & 2+1   & 237   &    & 3+1   & 206   &         & 3+1   & 824   &       & 3+1   & 818   &          & 3+1   & 216   &  \\
    \hline
    \multicolumn{3}{|c}{\footnotesize\# {solved (S)}} & \multicolumn{3}{|c}{16} & \multicolumn{3}{|c}{29} & \multicolumn{3}{|c}{29} & \multicolumn{3}{|c}{29} & \multicolumn{3}{|c}{29} & \multicolumn{3}{|c}{29} & \multicolumn{3}{|c}{29} & \multicolumn{3}{|c}{29} & \multicolumn{3}{|c}{29}  & \multicolumn{3}{|c|}{29} \\
    \hline
    \multicolumn{3}{|c}{\footnotesize{\# unsolved (U)}} & \multicolumn{3}{|c}{14} & \multicolumn{3}{|c}{1} & \multicolumn{3}{|c}{1} & \multicolumn{3}{|c}{1} & \multicolumn{3}{|c}{1} & \multicolumn{3}{|c}{1} & \multicolumn{3}{|c}{1} & \multicolumn{3}{|c}{1} & \multicolumn{3}{|c}{1}  & \multicolumn{3}{|c|}{1} \\
    \hline
    \multicolumn{3}{|c}{\footnotesize{avg. sec.: S}} & \multicolumn{3}{|c}{1734} & \multicolumn{3}{|c}{430} & \multicolumn{3}{|c}{246} & \multicolumn{3}{|c}{397} & \multicolumn{3}{|c|}{267} & \multicolumn{3}{|c}{275} & \multicolumn{3}{|c}{176}  & \multicolumn{3}{|c}{341} & \multicolumn{3}{|c|}{344}  & \multicolumn{3}{c|}{179} \\
    \hline
    \multicolumn{3}{|c}{\footnotesize{avg. gap: U}} & \multicolumn{3}{|c}{5.26} & \multicolumn{3}{|c}{37.24} & \multicolumn{3}{|c}{35.59} & \multicolumn{3}{|c}{33.29} & \multicolumn{3}{|c|}{38.01} & \multicolumn{3}{|c}{38.01} & \multicolumn{3}{|c}{37.03}  & \multicolumn{3}{|c}{37} & \multicolumn{3}{|c}{NA}  & \multicolumn{3}{|c|}{36.69} \\
    \hline
    \multicolumn{3}{|c}{\footnotesize {CPX/BD: XBS}}  & \multicolumn{3}{|c}{} & \multicolumn{3}{|c}{13.5} & \multicolumn{3}{|c}{18.4} & \multicolumn{3}{|c}{10.2} & \multicolumn{3}{|c}{14.8} & \multicolumn{3}{|c}{14.6} & \multicolumn{3}{|c}{22} & \multicolumn{3}{|c}{12.6} & \multicolumn{3}{|c}{13.3}   & \multicolumn{3}{|c|}{18.2} \\
    \hline
    \end{tabular}}%
  \label{base_bin}%
\end{table}%
\end{landscape}

\begin{landscape}
\begin{table}[htbp]
  \centering
  \caption{Comparison of Algorithm Performance under $\varepsilon = 0.1$ on Instances with General Integer $\mathbf{x}$}
    \scalebox{0.72}{\setlength{\tabcolsep}{0.01em} \begin{tabular}{|c|c|c|lcr|lcr|lcr|lcr|lcr|lcr|lcr|lcr|lcr|lcr|}
    \addlinespace
     \hline
    \multirow{2}[0]{*}{$K$} & \multirow{2}[0]{*}{} &  \multirow{2}[0]{-5pt}{} & \multicolumn{3}{|c|}{CPX} & \multicolumn{3}{c|}{BD0} & \multicolumn{3}{c|}{BD1} & \multicolumn{3}{c|}{BD2} & \multicolumn{3}{c|}{BD3} & \multicolumn{3}{c|}{BD4} & \multicolumn{3}{c|}{BD5} & \multicolumn{3}{c|}{BD6} & \multicolumn{3}{c|}{BD7} & \multicolumn{3}{c|}{BD8} \\
  \cline{4-33}
          &       &   & obj.   & sec.   & g(\%) & itr. & sec.  & g(\%)  & itr. & sec.  & g(\%) & itr. & sec.  & g(\%) & itr. & sec.  & g(\%) & itr. & sec.  & g(\%) & itr. & sec.  & g(\%) & itr. & sec.  & g(\%) & itr. & sec.  & g(\%) & itr. & sec.  & g(\%) \\
    \hline
    \multirow{10}[0]{*}{250} & \multirow{5}[0]{*}{$T_1$} & 1  & 1587.3   & 2050   &    & 3     & 93    &     & 3     & 64    &       & 3     & 89    &   & 5+1   & 50    &         & 5+1   & 51    &    & 3+1   & 35    &         & 3+1   & 42    &       & 3+1   & 44    &          & 3+1   & 60    &  \\
          &       & 2     & 1146.9 & T      & 14.72 & 9     & T     & 0.54  & 6     & 500   &       & 6     & 627   &    & 6+4   & 195   &        & 6+4   & 218   &    & 7+2   & 1602  &         & 7+3   & M & 0.52  & 7+2   & 1222  &          & 7+2   & M & 0.74 \\
          &       & 3     & 1652.4 & 3226   &    & 3     & M & 13.33    & 4     & 139   &       & 4     & 170   &    & 4+3   & 121   &        & 4+3   & 134   &   & 4+1   & 49    &          & 4+1   & 60    &       & 4+1   & 58    &          & 4+1   & 88    &  \\
          &       & 4     & 1683.2 & 770    &   & 2     & M & 15.46     & 3     & 1084  &       & 3     & 1093  &   & 4+1   & 49    &         & 4+1   & 49  & & 4+1   & 62    &             & 4+1   & 51    &       & 4+1   & 55    &          & 3+1   & 70    &  \\
          &       & 5     & 2387.8 & 559    &     & 4     & 88    &      & 3     & 100   &       & 3     & 120   &    & 5+1   & 45    &        & 5+1   & 45    &    & 4+1   & 37    &         & 4+1   & 54    &       & 4+1   & 53    &          & 3+1   & 78    &  \\
          \cline{2-3}
          & \multirow{5}[0]{*}{$T_2$} & 1  & 2832.3   & 1501   &    & 3     & 403   &      & 2     & 244   &       & 2     & 247   &    & 2+1   & 149   &         & 2+1   & 153   &    & 3+1   & 157   &        & 3+1   & 931   &       & 3+1   & 944   &          & 2+1   & 849   &  \\
          &       & 2     & 3087 & T      & 6.7  & 2     & M & 5.49   & 2     & M & 2.44  & 2     & M & 2.44  & 4+2   & 326   &     & 4+2   & 335   &   & 4+1   & 597   &         & 4+1   & M & 0.5   & 4+1   & M & 0.51      & 7+1   & M & 0.6 \\
          &       & 3     & 2864 & 1203   &   & 3     & 528   &       & 2     & 411   &       & 2     & 424   &   & 3+1   & 97    &         & 3+1   & 98    &   & 3+1   & 347   &          & 3+1   & 378   &       & 3+1   & 359   &          & 2+1   & 420   &  \\
          &       & 4     & 1853.6 & T      & 6.79 & 2     & M & 14.94 & 2     & M & 10.21 & 2     & M & 10.21 & 4+1   & M & 7.09 & 4+1   & M & 1.23  & 3+1   & M & 1.14  & 3+1   & M & 1.24  & 3+1   & M & 1.1     & 3+1   & M & 1.14 \\
          &       & 5     & 2173.5 & T      & 4.82   & 3     & M & 11.36 & 2     & 95    &       & 2     & 106   &   & 3+1   & 108   &         & 3+1   & 109   &   & 3+1   & 109   &         & 3+1   & 187   &       & 3+1   & 221   &           & 3+1   & 270   &  \\
          \cline{1-3}
    \multirow{10}[0]{*}{500} & \multirow{5}[0]{*}{$T_1$} & 1   & 3420.2  & T      & 14.58  & 6     & T  & 6.19  & 4     & M & 0.71  & 5     & T     & 0.7 & 3+1   & M & 0.69  & 3+1   & M & 0.69 & 3+4   & T     & 0.64 & 3+2   & M & 0.75  & 3+1   & M & 0.73      & 3+1   & M & 0.75 \\
          &       & 2     & 1566.5 & 1913   &    & 3     & 133   &       & 2     & 72    &       & 2     & 94    &   & 6+1   & 206   &          & 6+1   & 207   &    & 3+1   & 85    &          & 3+1   & 125   &       & 3+1   & 125   &        & 3+1   & 206   &  \\
          &       & 3     & 1209.2 & T      & 6.45  & 4     & M & 2.89   & 4     & M & 11.88 & 5     & M & 0.85  & 3+4   & 444   &    & 3+4   & 511   &   & 6+3   & M & 0.73       & 6+1   & M & 0.82  & 6+1   & M & 0.65    & 6+1   & M & 0.53 \\
          &       & 4     & 1511.7  & T     & 10.2  & 2     & M& 9.82  & 2     & M & 6.37  & 2     & M & 6.37  & 5+2   & 242   &       & 5+2   & 294   &    & 3+1   & 78    &         & 3+1   & 137   &       & 3+1   & 128   &        & 3+1   & 241   &  \\
          &       & 5     & 1013  & T      & 10.98   & 3     & M & 23.3  & 3     & M & 16.99 & 5     & T     & 0.7   & 3+3   & 330   &          & 3+3   & 387   &   & 5+3   & M & 0.65       & 5+1   & M & 0.74  & 5+1   & M & 0.77  & 5+1   & M & 0.68 \\
          \cline{2-3}
          & \multirow{5}[0]{*}{$T_2$} & 1  & 1043.6   & 1255   &     & 5     & T     & 18.68     & 2     & 139   &       & 2     & 177   &   & 4+1   & 229   &         & 4+1   & 228   &   & 3+1   & 174   &          & 3+1   & 319   &       & 3+1   & 382   &          & 3+1   & 510   &  \\
          &       & 2   & 1739.6  & T      & 10.33 & 3     & 270   &    & 2     & 111   &       & 2     & 153   &   & 2+1   & 135   &         & 2+1   & 135   &   & 3+1   & 163   &          & 3+1   & 291   &       & 3+1   & 314   &          & 2+1   & 395   &  \\
          &       & 3     & 2118.5 & 2896   &   & 5     & T     & 11.14     & 2     & 247   &       & 2     & 286   &   & 2+1   & 79    &         & 2+1   & 80    &    & 3+1   & 85    &         & 3+1   & 373   &       & 3+1   & 424   &          & 2+1   & 435   &  \\
          &       & 4     & 9625 & 949    &    & 3     & M & 0.2     & 2     & 1389  &       & 2     & 1385  &    & 2+1   & 1021  &        & 2+1   & 1020  &    & 3+1   & 875   &        & 3+1   & 1607  &       & 3+1   & 1607  &           & 2+1   & 1365  &  \\
          &       & 5     & 2037.7  & T     & 12.6 & 3     & 496   &   & 2     & 268   &       & 2     & 310   &   & 2+1   & 207   &          & 2+1   & 209   &   & 3+1   & 276   &         & 3+1   & 604   &       & 3+1   & 586   &          & 2+1   & 579   &  \\
         \cline{1-3}
    \multirow{10}[0]{*}{1000} & \multirow{5}[0]{*}{$T_1$} & 1   & 2499.4  & 3200   &    & 3     & 749   &       & 2     & 308   &       & 2     & 384   &  & 4+1   & 612   &          & 4+1   & 613   &   & 3+1   & 457   &        & 3+1   & 684   &       & 3+1   & 644   &            & 3+1   & 983   &  \\
          &       & 2   & 1400.2  & 2561   &   & 3     & M & 0.52    & 3     & 235   &       & 3     & 432   &    & 3+1   & 210   &        & 3+1   & 222   &    & 3+1   & 220   &         & 3+1   & 336   &       & 3+1   & 420   &          & 2+1   & 451   &  \\
          &       & 3   & 1298.4  & T      & 8.44 & 6     & T     & 15.21 & 3     & 535   &       & 3     & 791   &    & 3+2   & 765   &          & 3+2   & 931   &   & 4+1   & 449   &        & 4+1   & 382   &       & 4+1   & 462   &          & 3+1   & 933   &  \\
          &       & 4   & 3315  & T       & 8.1  & 5     & T     & 16.67 & 2     & 334   &       & 2     & 405   &   & 4+1   & 600   &         & 4+1   & 605   &   & 3+1   & 444   &         & 3+1   & 624   &       & 3+1   & 588   &           & 3+1   & 917   &  \\
          &       & 5   & 1149  & 2511   &   & 5     & T     & 14.6    & 2     & 144   &       & 2     & 214   &   & 4+1   & 348   &        & 4+1   & 351   &   & 3+1   & 334   &          & 3+1   & 346   &       & 3+1   & 345   &           & 2+1   & 543   &  \\
         \cline{2-3}
          & \multirow{5}[0]{*}{$T_2$} & 1 & 2552.7    & 3502   &   & 3     & 2524  &     & 2     & 1646  &       & 2     & 1639  &    & 2+1   & 987   &        & 2+1   & 999   &    & 3+1   & 1008  &          & 3+1   & 1268  &       & 3+1   & 1583  &         & 2+1   & 2613  &  \\
          &       & 2  & 1547.9   & T      & 7.71  & 3     & 1023  &   & 2     & 758   &       & 2     & 953   &    & 2+1   & 510   &        & 2+1   & 504   &   & 3+1   & 458   &           & 3+1   & 1143  &       & 3+1   & 1195  &         & 2+1   & 1499  &  \\
          &       & 3  & 2858.1   & T      & 7.51  & 5     & T     & NA  & 2     & 727   &       & 2     & 720   &   & 2+1   & 1001  &          & 2+1   & 993   &    & 3+1   & 531   &        & 3+1   & 1670  &       & 3+1   & 1418  &          & 2+1   & 1100  &  \\
          &       & 4   & 1778.1  & T      & 9.89 & 5     & T     & NA & 3     & 1300  &       & 3     & 1799  &   & 2+2   & 1248  &         & 2+3   & 2687  &    & 3+2   & 1765  &        & 3+2   & 3228  &       & 3+3   & T     & 4.14      & 2+1   & 2522  &  \\
          &       & 5   & 2700.7   & T     & 10.27 & 3     & 1210  &   & 2     & 921   &       & 2     & 1039  &    & 2+1   & 580   &       & 2+1   & 580   &    & 3+1   & 509   &        & 3+1   & 1690  &       & 3+1   & 1819  &            & 2+1   & 2032  &  \\
    \hline
          \multicolumn{3}{|c|}{\footnotesize{\# solved (S)}} & \multicolumn{3}{c|}{14} & \multicolumn{3}{c|}{11}  & \multicolumn{3}{c|}{24} & \multicolumn{3}{c|}{24} & \multicolumn{3}{c|}{28} & \multicolumn{3}{c|}{28} & \multicolumn{3}{c|}{26} & \multicolumn{3}{c|}{24} & \multicolumn{3}{c|}{24}  & \multicolumn{3}{c|}{24} \\
    \hline
    \multicolumn{3}{|c|}{\footnotesize{\# unsolved (U)}} & \multicolumn{3}{c|}{16} & \multicolumn{3}{c|}{19} & \multicolumn{3}{c|}{6} & \multicolumn{3}{c|}{6} & \multicolumn{3}{c|}{2} & \multicolumn{3}{c|}{2} & \multicolumn{3}{c|}{4} & \multicolumn{3}{c|}{6} & \multicolumn{3}{c|}{6}   & \multicolumn{3}{c|}{6} \\
    \hline
    \multicolumn{3}{|c|}{\footnotesize{avg. sec.: S}} & \multicolumn{3}{c|}{2007} & \multicolumn{3}{c|}{683}  & \multicolumn{3}{c|}{490} & \multicolumn{3}{c|}{569} & \multicolumn{3}{c|}{389} & \multicolumn{3}{c|}{455} & \multicolumn{3}{c|}{419} & \multicolumn{3}{c|}{689} & \multicolumn{3}{c|}{625}   & \multicolumn{3}{c|}{798} \\
    \hline
    \multicolumn{3}{|c|}{\footnotesize{avg. gap: U}} & \multicolumn{3}{c|}{9.38} & \multicolumn{3}{c|}{10.61} & \multicolumn{3}{c|}{8.1} & \multicolumn{3}{c|}{3.55}  & \multicolumn{3}{c|}{3.89} & \multicolumn{3}{c|}{0.96} & \multicolumn{3}{c|}{0.79} & \multicolumn{3}{c|}{0.76} & \multicolumn{3}{c|}{1.32}   & \multicolumn{3}{c|}{0.74} \\
    \hline
    \multicolumn{3}{|c|}{\footnotesize{CPX/BD: XBS}} & \multicolumn{3}{c|}{} & \multicolumn{3}{c|}{7.8} & \multicolumn{3}{c|}{11.4} & \multicolumn{3}{c|}{8.8} & \multicolumn{3}{c|}{14.2} & \multicolumn{3}{c|}{13.9} & \multicolumn{3}{c|}{18.5} & \multicolumn{3}{c|}{13.1} & \multicolumn{3}{c|}{12.8}   & \multicolumn{3}{c|}{9.1} \\
    \hline
    \end{tabular}}%
  \label{base_int}%
\end{table}%
\end{landscape}

\begin{table}[htbp]
  \centering
  \caption{Comparison of Algorithm Performance under $\varepsilon = 0.05$}
    \scalebox{0.72}{\begin{tabular}{|c|c|c|lcr|lcr|lcr|lcr|}
    \addlinespace
    \hline
    \multirow{3}[0]{*}{$K$} & \multirow{3}[0]{*}{} & \multirow{3}[0]{*}{} & \multicolumn{6}{|c|}{binary $\mathbf{x}$}                  & \multicolumn{6}{|c|}{integer $\mathbf{x}$}                 \\
    \cline{4-15}
          &       &       & \multicolumn{3}{|c|}{CPX} & \multicolumn{3}{|c|}{BD5} & \multicolumn{3}{|c|}{CPX} & \multicolumn{3}{|c|}{BD3} \\
    \cline{4-15}
          &       &       & obj.  & sec.  & g(\%) & itr.  & sec.  & g(\%) & obj.  & sec.  & \multicolumn{1}{c|}{g(\%)} & itr.  & sec.  & g(\%) \\
    \hline
    \multirow{10}[0]{*}{250} & \multirow{5}[0]{*}{$T_1$} & 1     & 2184.4 & 1084  &       & 4+1   & 39    &       & 1602.4 & 193   & \multicolumn{1}{c|}{} & 5+1   & 38    &  \\
          &       & 2     & 2053.8 & 68    &       & 3+1   & 22    &       & 1160 & T     & \multicolumn{1}{c|}{11.14 } & 6+4   & 173   &  \\
          &       & 3     & 1582.2 & 1346  &       & 6+1   & 59    &       & 1779.2 & 2813  &       & 4+4   & 127   &  \\
          &       & 4     & 1798.5 & 2196  &       & 6+1   & 71    &       & 1698.2 & 107   & \multicolumn{1}{c|}{} & 4+2   & 53    &  \\
          &       & 5     & 1258 & 80    &       & 3+1   & 23    &       & 2402.6 & 68    & \multicolumn{1}{c|}{} & 5+1   & 42    &  \\
    \cline{2-3}
          & \multirow{5}[0]{*}{$T_2$} & 1     & 1546.3 & T     & 0.9   & 5+1   & 90    &       & 2845.1 & 1487  & \multicolumn{1}{c|}{} & 2+1   & 157   &  \\
          &       & 2     & 3076.4 & T     & 17.55 & 6+2   & 1249  &       & 3140.4 & T     & \multicolumn{1}{c|}{16.74 } & 4+2   & 322   &  \\
          &       & 3     & inf & 1318  & inf & 2+1   & 71     & inf & 2879.3 & 1236  & \multicolumn{1}{c|}{} & 3+1   & 80    &  \\
          &       & 4     & inf & 2467  & inf & 3+1   & 58     & inf & 1882.5 & T     & \multicolumn{1}{c|}{4.51 } & 4+2   & 156   &  \\
          &       & 5     & 1782 & T     & 3.33  & 3+1   & 69    &       & 2186 & 2890  & \multicolumn{1}{c|}{} & 3+1   & 49    &  \\
     \cline{1-3}
    \multirow{10}[0]{*}{500} & \multirow{5}[0]{*}{$T_1$} & 1     & 1370.8 & T     & 21.61 & 8+1   & 327   &       & 3428.4 & 3427  & \multicolumn{1}{c|}{} & 3+1   & 147   &  \\
          &       & 2     & 1222.5 & T     & 18.86 & 3+1   & 70    &       & 1677.9 & T     & \multicolumn{1}{c|}{6.13 } & 6+2   & 254   &  \\
          &       & 3     & 918.6 & T     & 10.9  & 3+1   & 55    &       & 1394.7 & T     & \multicolumn{1}{c|}{17.98 } & 3+2   & 168   &  \\
          &       & 4     & 1918.8 & T     & 7.63  & 4+1   & 118   &       & 1508.4 & T     & \multicolumn{1}{c|}{5.15 } & 5+3   & 356   &  \\
          &       & 5     & 470.9 & T     & 15.55 & 3+1   & 57    &       & 1150.8 & T     & \multicolumn{1}{c|}{20.19 } & 3+3   & 247   &  \\
      \cline{2-3}
          & \multirow{5}[0]{*}{$T_2$} & 1     & 1478.2 & T     & 8.75  & 3+1   & 89    &       & 1063.7 & T     & \multicolumn{1}{c|}{1.65 } & 4+1   & 171   &  \\
          &       & 2     & 2288.6 & 243   &       & 3+1   & 94    &       & 1746.0 & 1355  & \multicolumn{1}{c|}{} & 2+1   & 120   &  \\
          &       & 3     & 1749.6 & T     & 0.51  & 3+1   & 101   &       & 2124.3 & 1284  & \multicolumn{1}{c|}{} & 2+1   & 109   &  \\
          &       & 4     & 1313.2 & T     & 0.57  & 3+1   & 66    &       & 9635.6 & 820   &       & 2+1   & 752   &  \\
          &       & 5     & 2458.5 & 1408  &       & 3+1   & 79    &       & 2044.4 & 3490  & \multicolumn{1}{c|}{} & 2+1   & 159   &  \\
     \cline{1-3}
    \multirow{10}[0]{*}{1000} & \multirow{5}[0]{*}{$T_1$} & 1     & 1686.2 & 1831  &       & 3+1   & 229   &       & 2512.9 & T     & \multicolumn{1}{c|}{5.45 } & 4+1   & 443   &  \\
          &       & 2     & 1423.8 & 1438  &       & 3+1   & 195   &       & 1411.9 & 1078  & \multicolumn{1}{c|}{} & 3+1   & 289   &  \\
          &       & 3     & 1218.7 & 1692  &       & 3+1   & 180   &       & 1702.4 & T     & \multicolumn{1}{c|}{29.32 } & 3+2   & 685   &  \\
          &       & 4     & 1479.7 & 919   &       & 3+1   & 108   &       & 3225.1 & T     & \multicolumn{1}{c|}{5.01 } & 4+2   & 676   &  \\
          &       & 5     & 1409.4 & T     & 31.29 & 3+1   & 215   &       & 1157.7 & 847   & \multicolumn{1}{c|}{} & 4+1   & 314   &  \\
      \cline{2-3}
          & \multirow{5}[0]{*}{$T_2$} & 1     & 1447 & 1852  &       & 3+1   & 212   &       & 2559.3 & 3454  & \multicolumn{1}{c|}{} & 2+1   & 870   &  \\
          &       & 2     & 2147.7 & 1548  &       & 3+1   & 262   &       & 1555.3 & 496   & \multicolumn{1}{c|}{} & 2+1   & 414   &  \\
          &       & 3     & 2342.7 & T     & 36.88 & 4+3   & T     & 27.62 & 2865  & T     & \multicolumn{1}{c|}{9.83 } & 2+1   & 506   &  \\
          &       & 4     & 2576.8 & 3538  &       & 3+3   & 2724  &       & 2375.6 & T     & \multicolumn{1}{c|}{34.24 } & 2+1   & 1082  &  \\
          &       & 5     & 2334  & 1053  &       & 3+1   & 203   &       & 2693.6 & 2748  & \multicolumn{1}{c|}{} & 2+1   & 491   &  \\
   \hline
    \multicolumn{3}{|c|}{\footnotesize{\# solved (S)}} & \multicolumn{3}{|c|}{17} & \multicolumn{3}{|c|}{29} & \multicolumn{3}{|c|}{17} & \multicolumn{3}{|c|}{30} \\
   \hline
    \multicolumn{3}{|c|}{\footnotesize{\# unsolved (U)}} & \multicolumn{3}{|c|}{13} & \multicolumn{3}{|c|}{1} & \multicolumn{3}{|c|}{13} & \multicolumn{3}{|c|}{0} \\
    \hline
    \multicolumn{3}{|c|}{\footnotesize{avg. sec.:S}} & \multicolumn{3}{|c|}{1417} & \multicolumn{3}{|c|}{246} & \multicolumn{3}{|c|}{1635} & \multicolumn{3}{|c|}{315} \\
    \hline
    \multicolumn{3}{|c|}{\footnotesize{avg. gap: U}} & \multicolumn{3}{|c|}{13.41} & \multicolumn{3}{|c|}{27.62} & \multicolumn{3}{|c|}{12.87} & \multicolumn{3}{|c|}{NA} \\
    \hline
    \multicolumn{3}{|c|}{\footnotesize{CPX/BD: XBS}} & \multicolumn{3}{|c|}{}  & \multicolumn{3}{|c|}{13.2} & \multicolumn{3}{|c|}{}  & \multicolumn{3}{|c|}{11.8} \\
    \hline
    \end{tabular}}%
  \label{small_e}%
\end{table}%

\begin{table}[htbp]
  \centering
  \caption{Computational Performances on Difficult Instances under $\varepsilon = 0.1$}
    \scalebox{0.73}{\begin{tabular}{|c|c|c|cc|lcr|lcr|lcr|}
    \addlinespace
    \hline
    \multirow{2}[0]{*}{$K$} & \multirow{2}[0]{*}{} & \multirow{2}[0]{*}{} & \multicolumn{2}{c}{CPX} & \multicolumn{3}{|c|}{BD1} & \multicolumn{3}{|c|}{BD3} & \multicolumn{3}{|c|}{BD5} \\
    \cline{4-14}
          &       &       & obj.  & g(\%)    & itr.  & sec.  & g(\%) & itr.  & sec.  & g(\%) & itr.  & sec.  & g(\%)\\
    \hline
    \multirow{10}[0]{*}{250} & \multirow{5}[0]{*}{$T_1$} & 1     & 3441.3 & 24.87  & 7     & 1778  &  & 6+4   & 963   &   & 7+1   & 466   &  \\
          &       & 2     & 1387.6 & 16  & 4     & 1192  &  & 5+3   & 213   &  & 5+1   & 227   &    \\
          &       & 3     & 3477.8 & 30.01  & 6     & 1928  &  & 6+3   & 1910  &  & 7+1   & 582   &  \\
          &       & 4     & 2892.7 & 43.16  & 4     & M     & 8.6  & 5+3   & 1213  &    & 5+1   & 296   &   \\
          &       & 5     &3798.4  &44.58  & 6     & T     & 3.68 &7+4   &T   &9.46  & 6+3   & T     & 4.12  \\
    \cline{2-3}
          & \multirow{5}[0]{*}{$T_2$} & 1     & 3129  & 26.67  & 5     & T     & 8.88  & 5+4   & T     & 6.47  & 8+3   & T     & 8.55   \\
          &       & 2     & 4657 & 51.2  & 2     & M     & 1.26  & 3+1   & 533   &  & 3+1   & 408   &  \\
          &       & 3     & 4477.6 & 36.33 & 6     & T     & 7.75  & 6+4   & T     & 11.42  & 7+4   & T     & 7.46  \\
          &       & 4     & 3237  & 51.35 & 7     & T     & 8.55  & 6+4   & T     & 13.25  & 8+3   & T     & 10.3 \\
          &       & 5     & 3971.9 & 45.27 & 5     & T     & 22.01 & 10+3  & T     & 20.96  & 6+3   & T     & 16.42  \\
     \cline{1-3}
    \multirow{10}[0]{*}{500} & \multirow{5}[0]{*}{$T_1$} & 1     & 1465.3 & 51.3  & 7     & T     & 10.07 & 7+5   & T     & 26.98  & 8+3   & T     & 16.57 \\
          &       & 2     & 3717.5 & 30.44  & 6     & T     & 6.5  & 4+4   & T     & 9.91  & 6+3   & T     & 11.36 \\
          &       & 3     & 1585.7 & 47.6  & 7     & T     & 0.92  & 5+5   & T     & 23.31  & 9+3   & T     & 11.4  \\
       &       & 4     & 2953.4 & 41.41  & 6     & T     & 8.48  & 6+4   & T     & 25.8  & 8+3   & T     & 11.27 \\
          &       & 5     & 3584.8 & 44.92 & 5     & T     & 20.62 & 6+4   & T     & 21.75  & 5+3   & T     & 16.35  \\
      \cline{2-3}
          & \multirow{5}[0]{*}{$T_2$} & 1     & 2190.3 & 33.5  & 2     & M     & 16.72  & 5+1   & 746   &  & 4+1   & M     & 0.93  \\
          &       & 2     & 3634.2 & 29.85   & 5     & T     & 12.32  & 6+4   & T     & 19.29   & 4+3   & T     & 7.41 \\
          &       & 3     & 3468.6 & 22.52  & 5     & T     & 0.55  & 4+1   &M & 8.6  & 4+3   & M     & 0.59  \\
          &       & 4     & 4310  & 28.54  & 3     & M     & 0.78  & 2+1   & M     & 0.53  & 3+1   & M     & 0.54  \\
          &       & 5     & 4314.74 & 64.22 & 5     & T     & 24.99 & 8+3   & T     & 25.63  & 4+3   & T     & 17.25  \\
     \cline{1-3}
    \multirow{10}[0]{*}{1000} & \multirow{5}[0]{*}{$T_1$} & 1     & 2117.7 & 23.04   & 6     & T     & 1.1  & 5+4   & T     & 9.78  & 5+3   & T     & 5.57 \\
          &       & 2     & 3273.5 & 53.6   & 5     & T     & 23.85  & 5+3   & T     & 30.14  & 4+3   & T     & 18.85  \\
          &       & 3     & 1806.2 & 64.91  & 5     & T     & 21.99 & 4+4   & T     & 46.79  & 5+3   & T     & 24.64 \\
          &       & 4    & 3679.1 & 48.57 & 5     & T     & 27.43  & 6+3   & T     & 35.99  & 4+3   & T     & 22.89  \\
          &       & 5     & 3866.25 & 52.92  & 5     & T     & 24.83  & 6+3   & T     & 35.03  & 4+3   & T     & 22.27 \\
      \cline{2-3}
          & \multirow{5}[0]{*}{$T_2$} & 1     & 1787.2 & 55.51   & 4     & T     & NA & 6+2   & M & NA & 5+3   & T     & NA \\
          &       & 2     & 3750.5 & 17.06   & 4     & T     & 10.2 & 3+3   & T     & 13.5  & 3+3   & T     & 10.51 \\
          &       & 3     & 2848.3 & 26.21    & 4     & 3205  &   & 5+3   & T     & 24.08  & 5+3   & T     & 9.89  \\
          &       & 4     & 3892.77 & 55.55  & 4     & T     & 33.97  & 5+1   & M     & NA  & 3+3   & T     & NA  \\
          &       & 5     & 6202.8 & 44.88  & 1     & M  & NA  & 3+1   & M     & 18.24 & 3+3   & T     & 15.61  \\
    \hline
    \multicolumn{3}{|c|}{\footnotesize{\# solved (S)}} & \multicolumn{2}{c}{0} & \multicolumn{3}{|c|}{4} & \multicolumn{3}{|c|}{6} & \multicolumn{3}{|c|}{5}\\
    \hline
    \multicolumn{3}{|c|}{\footnotesize{\# unsolved (U)}} & \multicolumn{2}{c}{30} & \multicolumn{3}{|c|}{26} & \multicolumn{3}{|c|}{24} & \multicolumn{3}{|c|}{25}\\
    \hline
    \multicolumn{3}{|c|}{\footnotesize{avg. sec.: S}} & \multicolumn{2}{c}{NA} & \multicolumn{3}{|c|}{2026} & \multicolumn{3}{|c|}{930} & \multicolumn{3}{|c|}{396}\\
    \hline
    \multicolumn{3}{|c|}{\footnotesize{avg. gap: U}} & \multicolumn{2}{c}{40.2} & \multicolumn{3}{|c|}{12.75} & \multicolumn{3}{|c|}{19.86} & \multicolumn{3}{|c|}{11.77}\\
    \hline
    \end{tabular}}%
  \label{hard}%
\end{table}%

\begin{table}[htbp]
  \centering
  \caption{Benders Decomposition and Jensen's Inequalities in OR Scheduling Problem}
    \scalebox{0.73}{\begin{tabular}{|c|c|c|ccc|ccc|ccc|ccc|}
    \addlinespace
    \hline
    \multirow{2}[0]{*}{Group} & \multirow{2}[0]{*}{} & \multirow{2}[0]{*}{} & \multicolumn{3}{c}{CPX} & \multicolumn{3}{|c|}{BD1} & \multicolumn{3}{c}{BD1J} & \multicolumn{3}{|c|}{BD1RJ} \\
    \cline{4-15}
          &       &       & obj.  & sec.  & g(\%) & itr.  & sec.  & g(\%) & itr.  & sec.  & g(\%) & itr.  & sec.  & g(\%) \\
    \hline
    \multirow{15}[0]{*}{I} & \multirow{5}[0]{*}{$\varepsilon=0$ (SP)} & 1     &20723.7	&20	&      & 171   & T     & 60.12  & 7     & 25    &       &       &       &  \\
          &       & 2     &20474.1	&23	&     & 166   & T     & 62.04  & 7     & 25    &       &       &       &  \\
          &       & 3     &20556.6	&22	&    & 162   & T     & 66.33  & 7     & 26    &       &       &       &  \\
          &       & 4     &20703.4	&25	&    & 161   & T     & 61.30  & 6     & 18    &       &       &       &  \\
          &       & 5     &20619.8	&18	&    & 161   & T     & 62.17  & 7     & 26    &       &       &       &  \\
          \cline{2-15}
          & \multirow{5}[0]{*}{$\varepsilon$ = 0.05} & 1     &19942.9	&T	&9.8  & 104   & T     & 68.37  & 8     & T     & 0.64   & 9     & T  & 0.57 \\
          &       & 2     &19730.1	&T	&11.64  & 104   & T     & 68.11  & 8     & T     & 0.66  & 9     & T  & 0.66 \\
          &       & 3     &19815.6	&T	&8.36  & 103   & T     & 68.16  & 9     & T & 0.69      & 8     & T  &  1.09\\
          &       & 4     &19937.4	&T	&21.06  & 103   & T     & 68.39  & 7    & T     & 0.7  &   9    & T & 0.67 \\
          &       & 5     &19873.5	&T	&29.88  & 102   & T     & 68.02  & 7     & T     & 0.8  & 8     & T     & 1.08\\
	   \cline{2-3}
          & \multirow{5}[0]{*}{$\varepsilon$ = 0.1} & 1     &19215.9	&T	&31.22  & 102   & T     & 68.09  & 6    & T     & 1.35  & 6     & T  & 3.77 \\
          &       & 2     &19049.7	&T	&15.19 & 103   & T     & 66.90  & 5     & T     & 1.48  & 6     & T &  4.5  \\
          &       & 3     &19120.2	&T	&38.04  & 99    & T     & 63.24  & 5     & T     & 1.28  & 6     & T  & 3.76 \\
          &       & 4     &19232.3	&T	&34.17  & 99    & T     & 67.24  & 6    & T     & 1.74  &  9     & T & 1.59 \\
          &       & 5     &19191.5	&T	&38.75  & 100   & T     & 66.22  & 5     & T     & 2.56  & 9     & T  & 2.76 \\
    \hline
    \multirow{15}[0]{*}{II} & \multirow{5}[0]{*}{$\varepsilon$ = 0 (SP)} & 1     &16566.9	&15	&    & 172   & T     & 49.50  & 5     & 15    &       &       &       &  \\
          &       & 2     &16433.5	&21	&     & 173   & T     & 54.44  & 5     & 15    &       &       &       &  \\
          &       & 3     &16551.8	&18	&    & 171   & T     & 54.39  & 4     & 11    &       &       &       &  \\
          &       & 4     &16507.3	&20	&  & 168   & T     & 48.77  & 5     & 15    &       &       &       &  \\
          &       & 5     &16256	&19	&     & 169   & T     & 57.38  & 5     & 20    &       &       &       &  \\
          \cline{2-15}
          & \multirow{5}[0]{*}{$\varepsilon$ = 0.05} & 1     &16091.7	&T	&5.86  & 106   & T     & 62.45  & 4     & M     & 1.9  & 6     & T  &  1.26\\
          &       & 2     &15942.4	&T	&5.51 & 108   & T     & 62.21  & 5    & T     & 2.0 & 6     & T  & 1.02 \\
          &       & 3    &16075.4	&T	&6.13  & 107   & T     & 65.73  & 6     & T     & 1.35  & 6     & T  &  1.14\\
          &       & 4     &16018.6	&T	&7.73  & 106   & T     & 59.25  & 5     & T     & 1.08  & 6     & T  & 1.1 \\
          &       & 5     &15823.5	&M	&5.2  & 109   & T     & 60.19  & 3    & M    & 6.86  & 6     & T  & 1.11 \\
          \cline{2-3}
          & \multirow{5}[0]{*}{$\varepsilon$ = 0.1} & 1     &15640.4	&T	&23.88  & 102   & T     & 61.31  & 5     & T     & 2.56  & 6     & T  & 2.65 \\
          &       & 2     &15496.3	&T	&17.87  & 97    & T     & 60.99  & 4     & T     & 6.95  & 6     & T  & 2.95 \\
          &       & 3     &15694.2	&M	&41.28  & 101   & T     & 63.47  & 5    & T     & 3.4  & 5     & T  &  2.47\\
          &       & 4    &15593.5	&T	&24.6  & 98    & T     & 55.18  & 4     & T     & 3.39  & 6 & T  &  2.27\\
          &       & 5     &15380.5	&T	&20.95  & 98    & T     & 58.86  & 4     & T     & 3.66  & 6     & T  &  2.43\\
      \hline
    \multicolumn{3}{|c|}{\footnotesize{avg. gap: U}} & \multicolumn{3}{|c|}{19.86} & \multicolumn{3}{|c|}{64.12} & \multicolumn{3}{|c|}{2.25} & \multicolumn{3}{|c|}{1.94} \\
    \hline
    \end{tabular}}%
  \label{JI}%
\end{table}%

%

\section{Conclusion}
\label{Sect_conclusion}

In this paper, we study chance constrained mixed integer program
with a finite discrete scenario set. We first present a
non-traditional bilinear formulation and analyze its structure. Its
linear counterparts are derived that could be stronger than popular
ones or big-M formulations. We also develop a variant of Jensen's
inequality that extends from stochastic program to chance
constrained program. To solve this type of challenging problem, we
provide a bilinear Benders reformulation and present a bilinear variant of Benders
decomposition method. 
In addition,  a few non-trivial enhancement strategies are designed to improve the solution capability of Benders
decomposition method. Different from existing
understanding that Benders decomposition might not be effective, we
observe that the presented implementation method,   jointly with appropriate enhancement techniques, performs
drastically better than a state-of-the-art commercial solver by an order of magnitude, on
deriving optimal solutions or reporting infeasibility.  Overall, the presented Benders decomposition
method basically does not depend on special assumptions and provides the first
easy-to-use fast algorithm to compute general chance constrained
program.

We believe that the presented bilinear formulation for chance
constrained program is informative and its structure is worth of
further analysis and investigation. One research direction is to develop
strong linear formulations and valid inequalities based on this bilinear
representation, especially for those whose recourse problem is an MIP.
We also point out that due to its simple scheme, the
presented Benders decomposition method is able to integrate almost
all improvement strategies and ideas developed for classical Benders
decomposition and stochastic programs.  So, it is rather an
algorithmic framework that can be further extended. Actually, we believe that if one improvement strategy benefits Benders decomposition method to solve SP, it would benefit to compute chance constrained program as well.  Clearly, more 
comprehensive study and evaluations are needed to support our understanding.



\bibliography{chance_constaint_refs}
\bibliographystyle{amsplain}

\newpage
\setcounter{equation}{0}
\renewcommand{\theequation}{A.\arabic{equation}}
\section*{Appendix}
\label{sect_appendix}
\noindent\textbf{Chance Constrained Operating Room Scheduling Problem}\\
The following chance constrained operating room scheduling problem is extended from its stochastic programming counterpart presented in \cite{batun2011operating}. The original parameters, notations and  descriptions are largely kept, while some additional variables and constraints are introduced to define the chance constrained model. Note that this formulation is  in \textbf{CC-bigM} form and  its bilinear \textbf{CC-MIBP}  form can be obtained straightforwardly.   


\vspace{2mm} \noindent \textbf{Indices}\\
\vspace{5mm}
\begin{tabular}{p{1.5cm}p{13cm}}
  $i,j$ & Surgery indices \\
  $k$ & Surgeon index\\
  $q,r$ & Operating rooperatingom indices  \\
  $\omega$ & Scenario index \\
  $i_k$  & Index of the first surgery of surgeon $k$\\
\end{tabular}

\vspace{2mm} \noindent \textbf{Parameters}\\
\vspace{5mm}
\begin{tabular}{p{1.5cm}p{13cm}}
  $L$ & Session length of each operating room\\
  $c^f$ & Fixed cost of opening an operating room\\
  $c^o$ & Overtime cost of an operating room per minute\\
  $c^S$ & Idle time cost of a surgeon per minute\\
  $s^S$ & Surgeon turnover time between two consecutive surgeries\\
  $s^R$ & Operating room trunover time between two consecutive surgeries\\
  $\pi_{\omega}$ & Probability of scenario $\omega$ \\
  $n$ & Total number of surgeries for scheduling\\
  $n_R$ & Total number of available operating rooms\\
  $n_S$ & Total number of surgeons\\
  $b_{ijk}$ & Binary parameter indicating whether surgery $i$ immediately precedes surgery $j$ in the surgery listing of surgeon $k$\\
  $pre_i(\omega)$ & Preincision duration of surgery $i$ in scenario $\omega$\\
  $p_i(\omega)$ & Incision duration of surgery $i$ in scenario $\omega$\\
  $post_i(\omega)$ & Postincision duration of surgery $i$ in scenario $\omega$\\
\end{tabular}

\vspace{2mm} \noindent \textbf{Decision Variables}\\
\vspace{5mm}
\begin{tabular}{p{1.5cm}p{13cm}}
  $x_{r}$ & Binary on/off status of operating room $r$ \\
  $y_{ir}$ & Binary variable indicating whether surgery $i$ is assigned to operating room $r$\\
  $z_{ijr}$ & Binary variable indicating whether surgery $i$ precedes $j$ in operating room $r$  \\
  $t_k$ & Continuous start time of surgeon $k$ \\
  $z^C(\omega)$  & Binary variable indicating whether scenario $\omega$ is not chosen as a responsive scenario\\
  $C_{ir}(\omega)$ & Continuous completion time of surgery $i$ in operating room $r$ under scenario $\omega$\\
  $I_{ij}(\omega)$ & Continuous idle time between surgery $i$ and $j$ in scenario $\omega$ (defined for $(i,j)$: $\sum_{k=1}^{n_S} b_{ijk} = 1$)\\
  $I_{k}(\omega)$ & Continuous idle time of surgeon $k$ before his/her first surgery in scenario $\omega$\\
  $O_r(\omega)$ & Continuous overtime in operating room $r$ in scenario $\omega$\\
  $I_{ij}^C(\omega)$  & Continuous auxiliary variable of $I_{ij}(\omega)$ for the chance constraint\\
  $I_{k}^C(\omega)$ & Continuous auxiliary variable of $I_{k}(\omega)$  for the chance constraint\\
  $O_r^C(\omega)$ & Continuous auxiliary variable of $O_r\omega)$  for the chance constraint \\
\end{tabular}
\newpage
\begin{align}
    \label{cobj}
    \displaystyle  \min \ & \sum_{r=1}^{n_R} c^f x_r + \sum_{\omega=1}^{K}\pi_{\omega}\mathcal{Q}^C(x,y,z,t,\xi(\omega))   \\
               \mbox{s.t.} \         &y_{ir} \leq x_r \qquad \forall i,r
                         \label{1b}\\
                         &\sum_{r=1}^{n_R} y_{ir} = 1 \qquad \forall i
                         \label{1c}\\
                         &z_{ijr} + z_{jir} \leq y_{ir} \qquad \forall i, j>i, r
                         \label{1d}\\
                         &z_{ijr} + z_{jir} \leq y_{jr} \qquad \forall i, j>i, r
                         \label{1e}\\
                         &z_{ijr} + z_{jir} \leq y_{ir}+y_{jr}-1 \qquad \forall i, j>i, r
                         \label{1f}\\
                         &t_k \leq L \qquad \forall k
                         \label{1g}\\
                         &x_{r} \geq x_{r+1} \qquad  r =1,...,n_R-1
                         \label{4a}\\
                         &\sum_{r=1}^{i} y_{ir} = 1 \qquad  i =1,...,min\{n, n_R\}
                         \label{4b}\\
                         &\sum_{q=r}^{min\{i, n_R\}} y_{iq} \leq \sum_{j=r-1}^{i-1} y_{j,r-1} \qquad \forall i, r \leq i
                         \label{4c}\\
                         &u_j \geq u_i + 1 - n(1-\sum_{r=1}^{n_R} z_{ijr}) \qquad \forall i, j \neq i
                         \label{5a}\\
                         &u_j \geq u_i + 1  \qquad \forall (i, j): \sum_{k=1}^{n_S} b_{ijk} = 1
                         \label{5b}\\
                         &\mathcal{Q}^C(x,y,z,t,\xi(\omega)) = \min \sum_{r=1}^{n_R} c^o O_r^C(\omega) + \sum_{(i,j): \sum_{k=1}^{n_S} b_{ijk} = 1} c^S I_{ij}^C(\omega) + \sum_{k=1}^{n_S} c^S I_k^C(\omega)
                         \label{3a}\\
                         &C_{ir}(\omega) \leq M y_{ir} \qquad \forall \omega, i, r 
                         \label{3b}\\
                         &C_{jr}(\omega) - C_{ir}(\omega) \geq s^R + pre_j(\omega) + p_j(\omega) + post_j(\omega) - M(1-z_{ijr}) \qquad \forall \omega,i,j \neq i,r 
                         \label{3c}\\
                         &I_k(\omega) - \sum_{r=1}^{n_R} C_{i_k r}(\omega) = -t_k - pre_{i_k}(\omega) - p_{i_k}(\omega) - post_{i_k}(\omega) \qquad \forall \omega,k 
                         \label{3d}\\
                         &\sum_{r=1}^{n_R} C_{ir}(\omega) \geq t_k + pre_{i}(\omega) + p_{i}(\omega) + post_{i}(\omega) \qquad \forall \omega,(i,k): \sum_{j=1}^{n} b_{jik}= 1 
                         \label{3e}\\
                         &I_{ij}(\omega) + \sum_{r=1}^{n_R} C_{ir}(\omega) - \sum_{r=1}^{n_R} C_{jr}(\omega) =  post_{i}(\omega) - s^S - p_j(\omega) - post_j(\omega) \notag \\
                         &  \   \qquad \qquad \qquad \qquad  \qquad   \qquad  \qquad  \qquad  \hspace{1cm} \forall \omega,(i,j): \sum_{k=1}^{n_S} b_{ijk}= 1 
                         \label{3f}\\
                         &O_r(\omega) - C_{ir}(\omega) \geq -L \qquad \forall \omega,i,r 
                         \label{3g}\\
                            &\sum_{\omega=1}^{K}\pi_{\omega}z^C(\omega) \leq \varepsilon
                         \label{ccz}\\ 
                         &I_{ij}^C(\omega) \geq  I_{ij}(\omega) - M z^C(\omega) \qquad \forall \omega,i,j
                         \label{cca}\\
                        &I_{k}^C(\omega) \geq  I_{k}(\omega) - M z^C(\omega) \qquad \forall \omega,k
                         \label{ccb}\\
                        &O_{r}^C(\omega) \geq  O_{r}(\omega) - M z^C(\omega) \qquad \forall \omega,r
                         \label{ccc}\\
\end{align}
\begin{align}
      %
                               &x_r \in \{0,1\}\  \forall r; y_{ir} \in \{0,1\}\  \forall i,r; z_{ijr} \in \{0,1\} \ \forall i, j \neq i, r; t_k \geq 0 \ \forall k \label{3h}  \\ 
                         &C_{ir}(\omega) \geq 0 \ \forall \omega,i,r; \ I_{ij}(\omega) \geq 0 \ \forall \omega,i,j; \ I_k(\omega)  \geq 0 \ \forall \omega,k;
                            \ O_r(\omega) \geq 0 \ \forall \omega,r \notag \\
                         &I_{ij}^C(\omega) \geq 0 \ \forall \omega,i,j; \ I_k^C(\omega)  \geq 0 \ \forall \omega,k;
                            \ O_r^C(\omega) \geq 0 \ \forall \omega,r; z^C(\omega) \in \{0,1\} \ \forall \omega  \notag
%
\end{align}

The objective function \eqref{cobj} is to minimize the sum of the (first-stage) fixed cost of opening operating rooms and the (second-stage) expected overtime and idle time cost from responsive scenarios. 
Constraints (\ref{1b}-\ref{1g}) impose restrictions on the first-stage variables.  Constraints (\ref{4a}-\ref{4c}) are introduced to eliminate symmetric solutions. Constraints (\ref{5a}-\ref{5b}) guarantee the feasibility 
of the first stage decisions. Constraints (\ref{3b}-\ref{3g}) impose restrictions on the second-stage variables. 

The inequality (\ref{ccz}) is the chance  constraint that limits the selection of non-responsive scenarios by $\varepsilon$. Noting that all of the cost coefficients are non-negative, we introduce 
constraints (\ref{cca})-(\ref{ccc}) to ensure that when a scenario is a  non-responsive one, there is no cost incurred from overtime or idle time of this scenario. 

Readers are encouraged to refer to  \cite{batun2011operating} for the detailed description and discussions on background, formulation and management insights of this particular application.
\end{document}